\documentclass[12pt,a4paper]{amsart}
\usepackage{listings}
\usepackage{amsmath}
\usepackage{amssymb}
\usepackage{graphicx}
\usepackage{longtable}
\usepackage{chapterbib}
\usepackage[timestamp,first]{draftcopy}
\usepackage[first]{draftcopy}
\draftcopySetScale{100}
\draftcopyPageX{200}
\draftcopyPageY{100}
\draftcopyBottomX{200}
\draftcopyBottomY{100}
\draftcopySetScaleFactor{0.8}

\usepackage[dvips]{epsfig}
\begin{document}
\title{\bfseries On Regularization and the Stokes Multiplier in Stirling's Approximation}
\author{\bfseries Victor Kowalenko}
\address{Department of Mathematics and Statistics\\ 
The University of Melbourne\\ Victoria 3010, Australia.}
\vspace{0.5 cm}
\email{vkowa@unimelb.edu.au}

\date{\today}



\begin{abstract}
Contrary to a recent work on the exactification of Stirling's approximation for the logarithm of the gamma function \cite{kow14}, 
Paris \cite{par14} has claimed that: (1) there is no need for the concept of regularization when determining 
values of $\ln \Gamma(z)$ from its complete asymptotic expansion and (2) the definition for the Stokes multiplier in the 
subdominant part of the complete expansion, which is responsible for demonstrating that the Stokes phenomenon is 
discontinuous rather than a smooth transition, is not correct. Here it is shown that unbeknownst to him, Paris has begun his 
analysis with regularized quantities. Then it is shown that the Stokes multiplier is entirely consistent with the 
conventional view of the Stokes phenomenon as described in Ref.\ \cite{ber89}. Finally, it is pointed out that 
the smoothing view as proposed by Paris is unable to reproduce the results in Table\ 7 of Ref.\ \cite{kow14}, 
which not only confirm the jump discontinuous nature of the Stokes phenomenon, but also provide accurate/hyperasymptotic values 
of $\ln \Gamma(z)$ to 30 figures, regardless of the size of the variable or whether the truncation is optimal or not.

\vspace{1.0cm}

{\bfseries Keywords:} Asymptotic series, asymptotic form, Borel summation, Cauchy integral, complete asymptotic expansion, conditional convergence, 
conventional view, discontinuity, divergence, exactification, gamma function, Mellin-Barnes regularization, regularization, remainder, 
smoothing, Stokes line, Stokes multiplier, Stokes phenomenon, Stokes sector, Stirling's approximation
\vspace{1.0cm}

{\bfseries 2010 Mathematics Subject Classification:} 30B10, 30B30, 30E15, 30E20, 34E05, 34E15, 40A05, 40G10, 40G99, 41A60

\end{abstract}
\maketitle
\section{Introduction}
Paris \cite{par14} has made several controversial statements concerning part of Sec.\ 3 in Ref.\ \cite{kow14}, which deals with 
the exactification of Stirling's approximation for the first time ever since its discovery almost three centuries ago. Exactification 
means that a complete asymptotic expansion for a function/integral has undergone regularization so that it is able to provide exact 
values of the original function for all values, including arguments or phases, of the main power variable. Despite acknowledging 
that the calculations in Ref.\ \cite{kow14} are basically correct, the two most important comments made by Paris are: (1) there is no need 
for the concept of regularization when determining values of $\ln \Gamma(z)$ from its complete asymptotic expansion and (2) the definition 
of the Stokes multiplier appearing in the subdominant part of the complete asymptotic expansion is not the correct interpretation of 
the quantity as envisaged by Stokes \cite{sto04} and Dingle \cite{din73}. The purpose of this article is not only to counter these claims, 
but also to present other fundamental points that could not be included in Ref.\ \cite{kow14}. As a consequence, the reader will be able 
to develop a better understanding of the primary issues in hyperasymptotics in which extremely accurate values, often to 30 figures or more, 
are sought from asymptotic expansions for functions/integrals. 

\section{Background}
Before Paris's comments can be addressed, we need to present a brief summary of Ref.\ \cite{kow14} to enable the reader to obtain
an understanding of the issues. The work begins with a derivation of the complete form of Stirling's approximation for $\ln \Gamma(z)$, 
which ultimately appears as (12) in Ref.\ \cite{kow14} and is given by
\begin{eqnarray}
\ln \Gamma(z) - \Bigl( z- \frac{1}{2} \Bigr) \ln z +z - \frac{1}{2} \,\ln (2 \pi) \equiv 2 z \sum_{k=1}^{\infty}
\frac{(-1)^{k+1}} {(2k-1)} \; \frac{\Gamma(2k) \, \zeta(2k)}{(2\pi z)^{2k}}  \;. 
\label{one}\end{eqnarray}
Usually, the asymptotic series on the rhs of the above statement, which is denoted by $S(z)$ in Ref.\ \cite{kow14}, is neglected 
and hence, the remaining terms on the lhs become the most familiar form of Stirling's approximation. On the few occasions that the 
series on the rhs is included, only the first few terms are usually given, while the remaining terms are neglected. That is, the asymptotic
series when displayed is often truncated at a few terms. The remainder is often denoted by $ + \cdots$ or by the 
Landau gauge symbol, $O()$. When this occurs, it is known as a standard Poincar$\acute{\rm e}$ asymptotic expansion. However, 
in order to obtain exact values of $\ln \Gamma(z)$ for any value of $z$, one requires not only all the terms in the above series, but 
also all the exponentially subdominant exponential terms that are generally discarded in standard Poincar$\acute{\rm e}$ asymptotics. Such 
terms are said to lie beyond all orders of a standard asymptotic expansion \cite{seg91}. The analysis leading to (\ref{one}) was not 
intended to be overly original, but to ensure that there were no exponentially subdominant terms in Stirling's approximation so that 
exactification could be achieved. In short, the above statement is a complete asymptotic expansion, but as we shall see later, for only one 
specific sector in the complex plane.  
 
The above statement possesses an equivalence symbol instead of an equals sign. This is because the infinite series on the
rhs is divergent in certain sectors of the complex plane and conditionally convergent for the remaining sectors. In the latter
case the equivalence symbol can be replaced by the more restrictive equals sign, but this is unnecessary since much of 
the analysis in Ref.\ \cite{kow14} is carried out with the equivalence symbol. In those sectors where the rhs is divergent, 
the finite values obtained on the lhs cannot possibly equal the values obtained from the infinite series on the rhs.
Therefore, the lhs is equivalent to the rhs, but it also means that the divergence must be tamed. The origin of the divergence 
is due to an impropriety in the method used to obtain the asymptotic series in the first place. That is, all asymptotic methods
are improper as described in great detail in Refs.\ \cite{kow09} and \cite{kow001}. Moreover, the divergence results in an infinity 
in the remainder of the asymptotic series. As a consequence, for the rhs of the above statement or equivalence to agree with the 
finite lhs, the infinity in the remainder must be removed in an appropriate manner, enabling the series to become summable. This 
process is known as regularization. Therefore, provided there is no truncation or the remainder has not been neglected, all 
asymptotic series must be regularized in order to yield a finite value. If they are truncated, then one is dealing with an 
approximation, which depending upon the value of $|z|$ can become extremely accurate. It is, therefore, important for the reader 
to realise that the finite truncated part represents an approximation, whose (relative) accuracy varies substantially throughout the complex plane.    
However, it can never equal the finite value on the lhs. For other values the approximation can be highly inaccurate, which means that 
one has ventured outside the range of applicability of the expansion. Hence exactification means that there are no ranges of applicability. 
That is, large or small values of the power variable are no longer necessary. Moreover, exactification can only be achieved with complete 
asymptotic expansions whereby all component asymptotic series must undergo regularization. 

Once the above equivalence is derived, Ref.\ \cite{kow14} proceeds to regularize $S(z)$ in order to obtain exact values of
$\ln \Gamma(z)$ for all values (including arguments) of $z$. To accomplish this, the Dirichlet series form for the 
Riemann zeta function is introduced and the order of the summations in the resulting equivalence are interchanged. As a
consequence, the asymptotic series on the rhs or $S(z)$ can be expressed as an infinite sum of Type I terminants.
Terminants were first introduced by Dingle in his remarkable book on asymptotic expansions \cite{din73}
because he noticed that many of the late terms in the asymptotic expansions of the special functions of mathematical physics could 
be approximated by them. Basically, he surmised that by studying the behaviour of the two types of terminants, one could come
up with very accurate approximations for the special functions of mathematical physics. It is still an open question whether
this can be achieved mainly because although the late terms become more accurate as an approximation to the asymptotic
expansion, the truncated part may already be too large as discussed in Ref.\ \cite{kow98}. Nevertheless, if a Type I terminant
can be regularized, then it follows that the rhs of the above equivalence can be regularized by introducing the regularized value
for each terminant in $S(z)$. This is basically the approach taken in the exactification of Stirling's approximation. The 
only problem with this approach is that because an infinite sum is involved, the regularized value of $S(z)$ may also become infinite, 
which is why the analysis in Ref.\ \cite{kow14} represents a departure from the analysis of asymptotic expansions with only one terminant
as carried out in Ref.\ \cite{kow09}. Thus, the theory in the latter reference, which represents an advance on the concepts and methods 
in Ref.\ \cite{din73}, was applied to $S(z)$. Once this was accomplished, the resulting infinite series was regularized as described in
the following sections. This enabled exact values of $\ln \Gamma(z)$ to be calculated for the first time ever via its asymptotic forms. 

As described in Ref.\ \cite{kow09}, there are at present two methods that are used widely for regularizing divergent series:
(1) Borel summation and (2) Mellin-Barnes (MB) regularization. A third technique based on the Euler-Maclaurin summation formula is
described in Ref.\ \cite{kow11}, but it has yet to be studied in detail. Moreover, logarithmically divergent series need
to be treated with special care \cite{kow11b,kow11c}, although from Lemma\ 2.2 in Ref.\ \cite{kow14}, one can see that the logarithmic power 
series can be regularized by using the regularized value of the geometric series, which is also a crucial step in Borel summation. 
Since there is a multitude of techniques that can be used evaluate an integral or solve a differential equation, there is more than
likely to be a host of methods for regularizing divergent series or asymptotic expansions, although the exactification of Stirling's approximation 
can be achieved by Borel summation and MB regularization. Hence there was no need to consider other methods of regularization in Ref.\
\cite{kow14}, which may not be the case for other problems in asymptotics.

The first example of exactification of a complete asymptotic expansion occurred in 1993 when T. Taucher and I carried out an extensive 
numerical investigation of the complete asymptotic expansion for the exponential series $S_3(a) = \sum_{n=1}^{\infty} \exp(-an^3)$,
which was given by 
\begin{align}
S_3(a) & \equiv 2 \sum_{k=0}^{\infty} \frac{(-1) ^{k+1} a^{2k+1}}{(2 \pi)^{6k+4}}\; \frac{\Gamma(6k+4)}{\Gamma(2k+2)}\;\zeta(6k+4) 
+ \frac{2 \sqrt{\pi}}{\Gamma(\frac{1}{6})\Gamma(\frac{5}{6})} \sum_{n=1}^{\infty} \frac{e^{-\sqrt{2z}}}{(6\pi n a)^{1/4}}
\nonumber\\ 
& \times \;\;  \sum_{k=0}^{\infty} \frac{\Gamma(k+1/6)}{(4 \sqrt{z})^{k}} \frac{\Gamma(k+5/6)}{\Gamma(k+1/2)} \;
\cos \Bigl( \sqrt{2z}\, -\frac{\pi}{8}- \frac{3k\pi}{4} \Bigr) .
\label{two}\end{align}
In the above equivalence $z \!=\! (2 n \pi/3)^3 a^{-1}$, while $\zeta(s)$ represents the Riemann zeta function. This presentation differs from
the original derivation in that the equivalence symbol has replaced an equals sign. However, since the numerical investigation in Ref.\ \cite{kow95}
was concerned with positive real values of $a$, where both series are deemed to be conditionally convergent, one could use an equals
sign in the analysis. As we shall see later, the regularization techniques used in Refs.\ \cite{kow14} and \cite{kow95} can be applied to
conditionally convergent asymptotic series. An asymptotic series is defined here as a power series with zero radius of absolute convergence. 

As a comparison, in standard Poincar$\acute{\rm e}$ asymptotics the above statement is expressed as
\begin{eqnarray}
S_3(a) \sim -\frac{2 \zeta(4)}{(2\pi)^4} \;\frac{\Gamma(4)}{\Gamma(2)}\; a + \frac{2 \zeta(10)}{(2\pi)^{10}}\, \frac{\Gamma(10)}{\Gamma(4)} \;
a^3 + O(a^5) \;\;.
\label{three}\end{eqnarray}
Such an expansion is referred to as a small $a$-expansion truncated at the third order, although exactly how small $a$ has to be for the
expansion to be valid is unclear. The second series in Equivalence\ (\ref{two}) is regarded as being exponentially subdominant to the first series 
because the factor of $\exp(-\sqrt{2z})$ vanishes as $a \to 0$. These terms are said to lie beyond all orders of the first or dominant asymptotic 
series \cite{seg91}. However, they are required if one wishes to obtain the exact values of $S_3(a)$ from Equivalence (\ref{two}). 

Both series in Equivalence\ (\ref{two}) diverge when $\pi/2 <|{\rm arg}\, z|< \pi$, but not in the positive half of the complex plane. Despite
the fact that they are conditionally convergent there, the coefficients diverge rapidly for values of $a$ greater than unity, which makes
the standard asymptotic form given by Approximation (\ref{three}) an inaccurate approximation when it is truncated. Moreover, even though both 
asymptotic series do not possess an infinity in their remainder, they can still undergo regularization except that now there is no need to
remove an infinity. Initially, our investigation concentrated on the Borel summation of both series, but it was found that the forms obtained 
via this method were not amenable to fast computation. So, we devised an alternative method of regularization, which cast the remainders of both 
series in the form of MB integrals. This enabled us to exploit the rapid exponential decay of the gamma function in the integrands 
along the imaginary axis. Consequently, we were able to obtain exact values of both series for values of $a$ ranging from 0.01 to 10 to incredible 
accuracy, in some cases as high as 65 figures. All the results from the MB-regularization of Equivalence (\ref{two}) were eventually presented 
in Secs.\ 7 and 8 of Ref.\ \cite{kow95} together with a comparison of the first step of the hyperasymptotic approach presented in Refs.\ 
\cite{ber91} and \cite{ber91a}. The latter approach was found to be far more inferior, especially for the larger values of $a$. Despite this 
remarkable achievement in that for the first time ever the exact values of a function had been obtained from its complete asymptotic expansion, the 
study has either gone largely unnoticed or even worse, been badly misunderstood.

As discussed in the preface of Ref.\ \cite{kow95}, in order to extend the analysis to the complex plane, a complete understanding 
of the Stokes phenomenon was required, which entailed not only determining where the lines of discontinuity occurred, but also the size of the 
subdominant jumps in an asymptotic expansion. This is an ambitious undertaking to say the least. Nevertheless, a major accomplishment towards 
this goal occurred with the publication of Ref.\ \cite{kow09}, even though it concentrated only on Borel-summable asymptotic series such as
both types of generalized terminants. These results were verified by deriving the corresponding MB-regularized forms. Thus, with the advent of 
Ref.\ \cite{kow09} it was now possible to consider the exactification of Stirling's approximation, which as indicated previously, involves an 
infinite number of terminants. Consequently, Ref.\ \cite{kow14} represents the next stage in the development of improved methods for handling 
asymptotic expansions over the past two decades of activity in asymptotics beyond all orders or hyperasymptotics as it has become more
familiarly known.

\section{Regularization}   
Now that the background material has been presented, we can now turn our attention to the main comments made by 
Paris \cite{par14}. These are: (1) there is no need for regularization in the analysis of asymptotic expansions and (2) the Stokes 
multiplier that arises when the asymptotic series on the rhs of Equivalence\ (\ref{one}) has been extended to the complex plane has not been 
interpreted correctly. This is despite the fact that Ref.\ \cite{kow14} has been able to present actual values of $\ln \Gamma(z)$ to 30
decimal places from the various asymptotic forms in the vicinity of the Stokes lines. Nevertheless, in order to discuss the second issue, we need
to address the first issue.

In Ref.\ \cite{par14} Paris expresses the main asymptotic series in the Stirling approximation as 
\begin{eqnarray}
\Omega(z) = \sum_{k=1}^{N-1} \frac{B_{2k}}{2k(2k-1)} + R_N(z) \;\;,
\label{four}\end{eqnarray}
where the remainder is given by  
\begin{eqnarray}
R_N(z) = \sum_{k=1}^{\infty} \frac{1}{k} \Bigl( e^{2pi i kz} T_{\nu} (2 \pi i k z) -e^{-2 \pi i kz} T_{\nu} (-2\pi i kz) \Bigl) \;\;,
\label{five}\end{eqnarray}
$T_{\nu}(z)$ is defined in terms of the incomplete gamma function as
\begin{eqnarray}
T_{\nu}(z) = \frac{e^{\pi i \nu}}{2\pi i} \, \Gamma(\nu) \Gamma(1-\nu,z) \;\;,
\label{fivea}\end{eqnarray}
and $\nu= 2N-1$. As stated in Ref.\ \cite{par14}, this result was first derived in Paris and Wood \cite{par92} via an MB integral approach. 

In Ref.\ \cite{kow14} the remainder is derived by a totally different approach based on Borel summation. There it appears as Eq.\ (83), which is
given by
\begin{align}
R_N^{SS}(z) & =\frac{ \Gamma(2N-1)}{2 \pi i} \sum_{n=1}^{\infty} \frac{1}{n} \Bigl( 
e^{-2 \pi n z i} \, \Gamma(2-2N, -2\pi n z i) 
\nonumber\\
& - \;\; e^{2 \pi n z i} \, \Gamma(2-2N, 2 \pi n z i ) \Bigr) \; .
\label{fiveb}\end{align}
Although these results appear to be different, Eq.\ (\ref{fiveb}) can be inferred from Eqs.\ (4.3), (4.10) and (4.11) in Ref.\ \cite{par92}.
That is, Eqs.\ (\ref{five}) and (\ref{fiveb}) are identical to each other. It should be noted here that not one result in Ref.\ 
\cite{par92} was ever used in Ref.\ \cite{kow14} because from the outset Equivalence\ (\ref{one}) was accompanied by the $\sim$ symbol instead
of an equivalence symbol. Hence Ref.\ \cite{par92} employs standard Poincar$\acute{\rm e}$ asymptotics, which is contrary to 
the presentation in Ref.\ \cite{kow14}.  

Ref.\ \cite{par14} also begins with a major inconsistency between (1.1) and (1.3). The first statement is written as a complete 
asymptotic expansion accompanied by the $\sim$ symbol, thereby implying divergence. Then the finite terms are given by the 
following definition:
\begin{eqnarray}
\Omega(z) := \log \Gamma(z) -(z-1/2) \log z +z - (1/2) \log 2 \pi \;\;.
\label{six}\end{eqnarray}
That is, $\Omega(z)$ is defined to be equal to the terms on the rhs of (\ref{six}). However, from Eq.\ (\ref{four}) it is also equal to the 
asymptotic series in the asymptotic expansion of $\ln \Gamma(z)$. Yet $\ln \Gamma(z)$ is written in terms of the same terms but with the 
$\sim$ symbol. Since $\Omega(z)$ is now finite, it must have a finite remainder for any value of $N$. The sudden replacement of the $\sim$ 
symbol by an equals sign together with the introduction of the definition symbol seems to have led Paris to conclude that there is no need for 
regularization. In fact, all that has been accomplished is that the complexity of the situation has been obscured.

Where the need for regularization has been missed is by the introduction of the so-called terminant function $T(z)$ in Eq.\ (\ref{five}).
As explained in Ch.\ 21 of Ref.\ \cite{din73}, terminants are divergent asymptotic series whose coefficients are given by the gamma function. 
For example, Dingle shows by Borel summation that a Type I terminant can be expressed as
\begin{eqnarray}
\sum_{k=n}^{\infty} \frac{\Gamma(k+ \alpha+1)}{(-x)^k} = \frac{\Gamma(n+\alpha +1)}{(-x)^n } \; \Lambda_{n+\alpha}(x) \;\;, 
\quad  |{\rm arg}\,  x|< \pi \;\;,
\label{seven}\end{eqnarray}    
where 
\begin{eqnarray}
\Lambda_s(x) =\frac{1}{\Gamma(n+1)} \int_0^{\infty} \frac{t^{s}\, e^{-t}}{1+ t/x} \; dt \;\;,
\label{sevena}\end{eqnarray} 
The integral in Eq.\ (\ref{sevena}) can be written in terms of the incomplete gamma function by using No.\ 3.383(10) in Ref.\ \cite{gra94}.
Then one finds that the above result reduces to
\begin{eqnarray}
\sum_{k=n}^{\infty} \frac{\Gamma(k+ \alpha+1)}{(-x)^k} = (-1)^n x^{\alpha-1} e^{x} \, \Gamma(n+\alpha) \, 
\Gamma(1-n-\alpha,x) \;\;,
\label{eight}\end{eqnarray}
where, in addition to the condition on ${\rm arg}\, x$, $\Re \, (n+ \alpha) >0$.

Eq.\ (\ref{eight}) also gives the impression that one does not need regularization, which may have led Paris to claim the same for his terminant 
function. Unfortunately, this is not correct. When $|{\rm arg}\, x| > \pi/2$, all we have to do is is consider 
$\Re\, (-x) >> \Im\, (-x)$  for the terminant, in which case the positive real part of $(-x)^k$ dominates the negative value produced by the 
imaginary part of $(-x)^k$. Consequently, all the terms in the terminant will be positive and thus, the series will become divergent. That is, an 
infinity exists in the remainder. Yet the rhs of Eq.\ (\ref{eight}) is finite. Hence there has been a ``sleight of hand" in deriving Eq.\ 
(\ref{eight}) via Borel summation.

In essence, Borel summation consists of a few steps. The first is that the gamma function on the lhs is replaced by its integral representation
and then the order of the summation and integration are interchanged. When this is done, one arrives at
\begin{eqnarray}
\sum_{k=n}^{\infty} \frac{\Gamma(k+ \alpha+1)}{(-x)^k} = \Bigl( \frac{-1}{x}\Bigl)^n  \int_0^{\infty} 
e^{-t} \,t^{n+ \alpha} \sum_{k=0}^{\infty} \Bigl(- \frac{t}{x}\Bigr)^k\; dt  \;\;.
\label{nine}\end{eqnarray}
Therefore, the Type I terminant has been expressed in terms of the geometric series. If we replace the geometric series by $1/(1+t/x)$, then we obtain Eq.\
(\ref{eight}). However, a major problem arises because according to p.\ 19 of Ref.\ \cite{whi73}, the geometric series is absolutely convergent for 
$|t/x| \!<\!1$ and divergent for $|t/x| \!\geq \!1$. Since $t$ ranges from zero to infinity, this substitution cannot be made for all values
of $t$. In short, ``Eq." (\ref{eight}) is simply not valid.

The situation is even more complex because the geometric series is not always divergent outside its radius of absolute convergence. To see 
this more clearly, let us replace $-t/x$ in the geometric series by $z$ for the time being. Consequently, we can use the material in Ch.\ 4 of Ref.\ 
\cite{kow09} or in Ref.\ \cite{kow001}, where the series is expressed as 
\begin{eqnarray}
\sum_{k=0}^{\infty}z^k =\sum_{k=0}^{\infty} \Gamma(k+1)\, \frac{z^k}{k!}=\lim_{p \to \infty} 
\sum_{k=0}^{\infty} \frac{z^k}{k!} \int_0^p  e^{-t} \,t^k \; dt \;\;. 
\label{ten}\end{eqnarray}
Since the integral in Eq.\ (\ref{ten}) is finite, technically, we can interchange the order of the summation
and integration. In reality, an impropriety has occurred when we do this, which will be discussed in more detail
shortly. For now, interchanging the summation and integration yields
\begin{align}
\sum_{k=0}^{\infty}z^k = & \lim_{p \to \infty} \int_0^p  e^{-t} \sum_{k=0}^{\infty} \frac{(zt)^k}{k!} \; dt=
\lim_{p \to \infty} \int_0^p e^{-t(1-z)} \; dt
\nonumber\\
& =\;\; \lim_{p \to \infty}\left[ \frac{e^{-p(1-z)}}{z-1}+ \frac{1}{1-z} \right] . 
\label{eleven}\end{align}
For $\Re\, z \!<\! 1$, the first term in the last member of the above equation vanishes and the geometric series yields a
finite value of $1/(1-z)$. Therefore, we see that the same limit value is obtained for the series when $\Re \, z \!<\!1$
as when $z$ lies in the circle of absolute convergence given by $|z| \!<\!1$ or the unit disk. According to the definition
on p.\ 18 of Ref.\ \cite{whi73}, this means that the series is conditionally convergent for $\Re \, z \!<\! 1$
and $|z| \!>\! 1$. For $\Re\, z \!>\!1$, however, the first term in the last member of Eq.\ (\ref{eleven})
is infinite. Since regularization is the process of removing the infinity in the remainder of a divergent series,
we remove or neglect the first term of the last member of Eq.\ (\ref{eleven}). As a consequence, we are left with a
finite part that equals $1/(1-z)$, which is known as the regularized value. Therefore, for all complex values of $z$
except $\Re \, z \!=\!1$, we have
\begin{eqnarray}
\sum_{k=0}^{\infty} z^k 
\begin{cases}
\equiv 1/(1-z)\;\;, & \;\; \Re \,z >1\;\;, \\
= 1/(1-z) \;\;, & \;\; \Re \, z <1.
\end{cases}
\label{twelve}\end{eqnarray}

At the barrier of $\Re \,z \!=\! 1$, the situation appears to be unclear. For $z \!=\!1$ the last member of
Eq.\ (\ref{eleven}) vanishes, which is consistent with removing the infinity from $1/(1-z)$. For other
values of $\Re\, z \!=\!1$, the last member of Eq.\ (\ref{eleven}) is clearly undefined. This is to be
expected as this vertical line represents the border between the domains of convergence and divergence for the geometric series.
Nevertheless, because the finite value remains the same to the right and to the left of the barrier at $\Re \, z \!=\!1$
and in keeping with the fact that regularization is effectively the removal of the first term on the rhs
of Eq.\ (\ref{eleven}), we take $1/(1-z)$ to be the finite or regularized value when $\Re \,z \!=\!1$.
Hence Equivalence\ (\ref{twelve}) becomes
\begin{eqnarray}
\sum_{k=0}^{\infty} z^k 
\begin{cases}
\equiv 1/(1-z)\;\;, & \;\; \Re \,z \geq 1\;\;, \\
= 1/(1-z)\;\;, & \;\; \Re \, z <1.
\end{cases}
\label{fourteen}\end{eqnarray}

Frequently, the values of $z$ over which an asymptotic series is convergent and those over which it is divergent are not 
known. So, in these cases it is better to replace the equals sign by the less stringent equivalence symbol on
the understanding that we may be dealing with a series that is absolutely or conditionally convergent for
some values of $z$. As a result, the above statement can be succinctly expressed as 
\begin{eqnarray}
\sum_{k=N}^{\infty}z^k =z^N \sum_{k=0}^{\infty} z^k \equiv \frac{z^N}{1-z} \;\;.
\label{fifteen}\end{eqnarray}
Such a statement is no longer an equation, but an equivalence statement or equivalence for short. Moreover, it can be extended to the 
binomial series as explained in Ch.\ 4 of Ref.\ \cite{kow09}. It should also be noted that the 
above result is only applicable when the form for the regularized value of the divergent series is identical to the form 
of the limiting value of the absolutely convergent series. Later in Ch.\ 4 of Ref.\ \cite{kow09} the regularization of
the hypergeometric series $_2F_1(a+1,b+1;a+b+2-x;1)$ is discussed, which is found to possess a regularized value for $\Re\, x>0$ 
that is different from the limit when the series is convergent for $\Re\, x<0$. Thus, in this instance the two forms must be 
kept separate from each other as in (\ref{fourteen}).

Now let us see what happens when we replace $z$ by $-t/x$. According to the above analysis, the geometric series in Eq.\ (\ref{nine}) is
divergent whenever $\Re\, (-t/x) >1$. As $t$ ranges from zero to infinity there will always be values where it is divergent for 
$\Re\, x<0$. For these values of $t$ and $x$ we require the regularized value of the geometric series because the series
is outside the barrier of convergence. Rather than determine the exact regions where the convergent and divergent regions apply,
we simply introduce Equivalence\ (\ref{fifteen}) into Eq.\ (\ref{nine}). Then we obtain
\begin{eqnarray}
\sum_{k=n}^{\infty} \frac{\Gamma(k+ \alpha+1)}{(-x)^k} \equiv \Bigl( \frac{-1}{x}\Bigl)^n  \int_0^{\infty} 
\frac{t^{n+ \alpha}\, e^{-t}}{(1+t/x)} \; dt \;\;,
\label{sixteen}\end{eqnarray}   
for $|{\rm arg}\, x| \neq \pi$. In other words, we have obtained Dingle's result given by Eq.\ (\ref{seven}) here, but 
with one major difference. That is, we now have an equivalence symbol instead of an equals sign, which informs us that the series
may be divergent for certain values of $x$. Even when it is convergent, the series may only be conditionally convergent, which is
equally troublesome to treat. Therefore, we have seen that regularization is crucial for handling terminants.

It should also be noted that if the steps in Borel summation are reversed, then an asymptotic expansion is obtained. That is,
if we wish to obtain an asymptotic series from the convergent integral on the rhs of ``Eq."\ (\ref{eight}), then all we need
to do is expand the denominator as a power series. This method of deriving an asymptotic series is known as the method of
expanding most of the exponential and is discussed on p.\ 113 of Ref.\ \cite{din73}. Like all asymptotic methods
it has an impropriety, which arises from using a power series expansion outside its radius of absolute convergence. The divergence
in a series expansion obtained by employing a standard technique such as Laplace's method or the method of steepest descent is 
therefore an indication that something improper has occurred. Regularization represents the necessary corrective measure for rendering
the values of the original function. Hence this is why it is necessary that the complete Stirling approximation must be regularized 
to yield the exact values of $\ln \Gamma(z)$ in Ref.\ \cite{kow14}.

Since we have seen that terminants can become divergent, the question now becomes: Why does Paris seem to think that regularization
is not needed at all in Ref.\ \cite{par14}? The answer to this question is that he has not used Borel summation to obtain the Borel-summed 
forms given in Ref.\ \cite{kow14}. Rather he has used MB integrals in his analysis. In fact, unbeknownst to him he has employed MB-regularization,
which becomes evident when one examines his ``Eq".\ (2.1). Before we can do this, however, we need to understand what MB 
regularization is. So, let us consider the general series given by
\begin{eqnarray}
S_{I}(N,z) =\sum_{k=N}^{\infty} f(k) (-z)^k \;\;,
\label{sixteena}\end{eqnarray}
where $N$ is referred to as the truncation parameter.

Basically, MB regularization applies whenever the function $f(k)$ in Eq.\ (\ref{sixteena})  possesses the following properties: 
(1) as $ L \! \to \! \infty$, $|f(s)| \!=\!O\bigl( \exp(-\epsilon_1 L)\bigr)$ for $s \!=\!c \! +\! iL$, and $|f(s)| \!=\! 
O \bigl( \exp(-\epsilon_2 L) \bigr)$ for $s \!=\! c \!-\!i L$, where $\epsilon_1, \epsilon_2 \!>\!0$, (2) $-\pi \!<\! \theta \!=\! 
{\rm arg}\, z \!<\! \pi$, (3) there exists a real number $c$ such that the poles of $\Gamma(N-s)$ lie to the right of the line 
$N \!-\!1  \!<\! c =\Re\,s  \!<\! N$ in the complex plane and that the poles of $f(s) \Gamma(s+1-N)$ to the left of it and (4) $
z^s f(s) \Gamma(s+1-N) \Gamma(N-s)$ is single-valued to the right of the line. 

Consider the following integral along the imaginary axis given by  
\begin{eqnarray}
I=(-1)^N \int_{c-i\infty}^{c+i\infty}  z^s f(s) \Gamma(1+s-N) \Gamma(N-s)\; ds \;\;,
\label{seventeen}\end{eqnarray}
where the cut-off $c$ is given above. Eq.\ (\ref{seventeen}) is basically an MB integral passing through $c$ when the line contour 
intersects the real axis. The first two properties of $f(s)$ are required for ensuring that the modulus of the integrand in the above 
integral decays exponentially at the endpoints. That is,
\begin{eqnarray}
\left| \frac{z^s \, f(s)}{e^{-i \pi s}-e^{i \pi s}} \right|  \stackrel{s=c \pm i L}{\sim} \, |z|^c 
e^{\mp L\theta} e^{-\pi L}  \, |f(c \pm iL) | \;\;.
\label{eighteen}\end{eqnarray}
Here, the $\sim$ symbol denotes "goes as". Therefore, the upper limit given by $s \!=\! c \!+\! i \infty$ decays exponentially for all 
values of $\theta$ in the principal branch provided $f(c \!+\! iL)| \!= \! O(\exp(-\epsilon_1 L))$ as $ L \! \to \! \infty$ and $\epsilon_1 \!>\! 0$. 
Similarly, the lower limit decays exponentially provided $|f(c \!-\! iL)|=O(\exp(-\epsilon_2 L))$ as $L \! \to \!  
\infty$ and $\epsilon_2 \!>\!0$. Furthermore, the integrand of the above integral is single-valued because ${\rm arg}\, z$ has 
been confined to the principal branch of the complex plane from the second property. Since the integrand is single-valued 
to the right of the line ${\rm Re}\,s \!=\! c$, we can close the line contour integral to the right and apply Cauchy's residue theorem. 
Hence we arrive at
\begin{eqnarray}
\int_{c-i\infty}^{c+i\infty} z^s f(s) \Gamma(1+s-N) \Gamma(N-s) \; ds + \int_{C} z^s f(s) \Gamma(1+s-N) \Gamma(N-s)\; ds
\nonumber\\
= \; 2 \pi i \sum {\rm Res}\left\{z^s f(s) \Gamma(1+s-N) \Gamma(N-s)
\right\} \;\;,
\label{eighteena}\end{eqnarray}
where the contour $C$ represents the great arc contour integral closing the limits of the MB integral, while $\sum {\rm Res} \{f(s)\}$ 
denotes the sum of the residues of $f(s)$ in the region bounded by the line $\Re\,  s=c$ and the arc contour. 

In accordance with the third property of $f(s)$, only the poles of $\Gamma(N \!-\! s)$ lie in the region to the right of the line contour. Therefore, 
the residues in Eq.\ (\ref{eighteena}) are the simple poles of the gamma function for all positive integers greater than or equal to $N$. Hence 
Eq.\ (\ref{eighteena}) reduces to
\begin{align}
\int_{c-i\infty}^{c+i\infty} z^s f(s) \Gamma(1+s-N) \Gamma(N-s) \; ds & +\int_{C} z^s f(s) \Gamma(1+s-N) \Gamma(N-s) \; ds
\nonumber\\
& = \; \; 2 \pi i (-1)^N \sum_{k=N}^{\infty} (-z)^k f(k) \;\;.
\label{nineteen}\end{align}

It should be noted that for the case when the series on the rhs of Eq.\ (\ref{nineteen}) is convergent, we have two cases: either the 
MB integral on the lhs and the integral along the great arc are both convergent or they are both divergent. Whichever of the two cases 
applies depends on the magnitude of ${\rm arg}\,z$. For example, consider the MB integral for the exponential function, which is
\begin{eqnarray}
\exp(-z)=\sum_{k=0}^{\infty}\frac{(-z)^k}{k!}= \frac{1}{2\pi i} \int_{\scriptstyle{c-i\infty}\atop c= {\rm Re}\,s>0}^{c+i\infty}
\Gamma(s) z^{-s}\; ds \;\;.
\label{twenty}\end{eqnarray}
The above equation appears on p.\ 261 of Ref.\ \cite{cop76} except that $s$ has been replaced by $-s$. According to Copson, this result
is valid only for $|{\rm arg}\,z|<\pi/2$. When evaluating the MB integral numerically by using the NIntegrate routine in Mathematica \cite{wol92},
one obtains values of $\exp(-z)$ very quickly for $|{\rm arg}\,z| \!<\! \pi/4$, but as ${\rm arg}$ approaches $\pi/2$, the evaluation 
becomes very slow, indicating that convergence problems are arising. For $|{\rm arg}\,z| \!>\! \pi/2$, the MB integral no longer converges. 
This means that we have encountered the second case where the integral along the great arc is also divergent.  Hence we have seen that the 
integral along the great arc does not always vanish even when the series on the rhs converges.

For the case where the complex series is divergent, Eq.\ (\ref{twenty}) can be written as
\begin{eqnarray}
\sum_{k=N}^{\infty} (-z)^k f(k) -\frac{(-1)^N}{2\pi i}\int_{C} z^s f(s) \Gamma(1+s-N) \Gamma(N-s)\; ds
\nonumber\\
= \; \frac{(-1)^N}{2\pi i} \int_{c-i\infty \atop{N-1 < c={\rm Re}\, s<N}}^{c+i\infty}z^s f(s) \Gamma(1+s-N) \Gamma(N-s) \; ds \;\;.
\label{twentyone}\end{eqnarray}
If the MB integral on the rhs of Eq.\ (\ref{twentyone}) is defined, i.e., it yields a definite value, then the integral along 
the great arc must also be divergent or else the divergence in the series cannot be removed. By removing the integral in this
situation, we are effectively regularizing the series. As noted in Ref.\ \cite{kow95}, this is somewhat analogous to evaluating 
the Hadamard finite or regularized part of a divergent integral \cite{lig75}. So, let us consider the evaluation of the contour 
integral along the great arc, which can be expressed as
\begin{eqnarray}
I_{{\text arc}}=  \lim_{L \to \infty}i \pi L \int_{-\pi/2}^{\pi/2} d \gamma \; 
\frac{f(Le^{i\gamma})\,e^{i\gamma}}{\sin(\pi Le^{i\gamma})} \; \exp(L e^{i\gamma}(\ln|z|+i \theta))\;\;.
\label{twentytwo}\end{eqnarray}
In obtaining this result we have used the reflection formula for the gamma function, which is
\begin{eqnarray}
\Gamma(1+s-N) \Gamma(N-s)=\frac{(-1)^{N+1}\pi}{\sin(\pi s)}\;\;.
\label{twentythree}\end{eqnarray}
The magnitude of the integral in Eq.\ (\ref{twentytwo}) is bounded by
\begin{align}
|I_{{\text arc}}| \leq 2 \pi &  \lim_{L \to \infty}L  \int_{-\pi/2}^{\pi/2}d\gamma \; |f(Le^{i\gamma})|\, 
\exp \left(L \ln|z| \cos \gamma \right)
\nonumber\\
& \times \;\; \exp\left( -L \pi |\sin \gamma|-L\theta \sin \gamma \right) \;\; .
\label{twentyfour}\end{align}
The above integral can be split into two integrals, one between $(-\pi/2,0)$ and the other between $(0,\pi/2)$. Then the Jordan inequality
can be applied. At the lower limits of both integrals there is no contribution from the $\sin \gamma$ terms in the exponential 
leaving only the terms involving $\cos\gamma$, which are bounded by unity in the integrals. Therefore, the conditions under which the 
contour along the great arc vanishes are $\epsilon_1 \!>\! \ln |z|$ and $\epsilon_2 \!>\! \ln |z|$. For $\epsilon_1 \!=\! \ln|z|$ 
and $\epsilon_2 \!=\! \ln|z|$, the integral can vanish but this will depend upon the algebraic behaviour of $f(s)$. For 
$\epsilon_1 \!<\! \ln|z|$ or $\epsilon_2 \!<\! \ln|z|$, $I_{{\text arc}}$ will yield infinity. In this case we simply remove the infinity 
in accordance with the process of regularization. Therefore, irrespective of whether the contour integral along the great arc yields 
infinity or vanishes, we arrive at the following equivalence statement:
\begin{align}
S_{I}(N,z)=\sum_{k=N}^{\infty} & \; f(k) \, (-z)^k \equiv  \frac{(-1)^N}{2\pi i} \int_{\scriptstyle{c-i\infty} 
\atop{N-1 < c=\Re s<N}}^{c+i\infty} z^s f(s)
\nonumber\\
& \times \;\; \Gamma(1+s-N) \, \Gamma(N-s) \;ds \;\;.
\label{twentyfive}\end{align}
We can simplify the above result by introducing again the reflection formula for the gamma function or Eq.\ (\ref{twentythree}). 
As a consequence, we arrive at
\begin{eqnarray}
S_{I}(N,z) =\sum_{k=N}^{\infty} f(k) (-z)^k \equiv \int_{c-i \infty \atop{N-1<c=\Re\,s<N}}^{c+i \infty}  
\frac{z^{s} \,f(s)}{e^{-i \pi s}-e^{i \pi s}} \;ds \;\;.
\label{twentysix}\end{eqnarray}
The conditions on $z^sf(s)$ as $s \!\to\! c \pm i\infty$ are required to ensure that the integral in Eq.\ (\ref{twentysix}) is convergent. 
They have been included more for the sake of completeness rather than utility because the above result has to be integrated numerically in
order to obtain a definite value. If the integral on the rhs of Eq.\ (\ref{twentysix}) is divergent, then numerical integration will fail.

Now we let $f(k) \!=\! |B_{2k}|\,z/(2k(2k-1)$ and $z \!=\! 1/z^2$ in $S_I(N,z)$ as given by Eq.\ (\ref{sixteena}). In order to determine the MB integral 
involving this $f(k)$, we have to continue the Bernoulli numbers analytically to complex values. This is done by expressing the Bernoulli 
numbers in terms of integer values of the Riemann zeta function and analytically continuing the latter to complex values. From No.\ 9.616 
of Ref.\ \cite{gra94}, we have
\begin{eqnarray}
B_{2k} =\frac{(-1)^{k-1} \,(2k)!}{2^{2k-1}\,\pi^{2k}}\; \zeta(2k)\;\;.
\label{twentyseven}\end{eqnarray}
Hence $f(k)= \Gamma(2k-1) \zeta(2k)/2(2\pi)^{2k} z$. Then according to Eq.\ (\ref{twentyone}), we obtain  
\begin{align}
\sum_{k=N}^{\infty} \frac{B_{2k}}{2k(2k-1)z^{2k-1}} & = \frac{1}{2\pi i} \int_{c-i\infty \atop{ {\rm Max}[1/2, N-1]< c={\rm Re}\, s< N}}^{c+i\infty} 
\frac{1}{(2\pi z)^{2s-1}} \; \frac{\Gamma(2s-1) \zeta(2s)}{\sin(s \pi)} \; ds
\nonumber\\
& + \;\; \frac{1}{2\pi i} \int_{C} \frac{1}{(2\pi z)^{2s-1}} \; \frac{\Gamma(2s-1) \zeta(2s)}{\sin(s \pi)} \; ds \;\;.
\label{twentyeight}\end{align}
Note that the lower limit has been adjusted. This is required in order that the singularities of the integrand other than those for $\Gamma(N-s)$
lie to the left of the line contour. Specifically, the line contour needs to lie to the right of the pole at $s=1/2$ in $\Gamma(2s-1)$, while 
the truncation parameter $N$ must be greater than zero. We now make a change of variable by putting $s=(s^{'}+1)/2$ in the MB integral on the 
rhs of the above equation. This yields
\begin{align}
\sum_{k=N}^{\infty} \frac{B_{2k}}{2k(2k-1)z^{2k-1}} & = \frac{1}{4\pi i} \int_{c^{'}-i\infty}^{c^{'}+i\infty} \frac{1}{(2\pi z)^{s^{'}}} \; 
\frac{\Gamma(s^{'}) \zeta(s^{'}+1)}{\sin((s^{'}+1) \pi/2)} \; ds^{'}
\nonumber\\
& + \;\; \frac{1}{2\pi i} \int_{C} \frac{1}{(2\pi z)^{2s-1}} \; \frac{\Gamma(2s-1) \zeta(2s)}{\sin(s \pi)} \; ds \;\;,
\label{twentyeighta}\end{align}
where ${\rm Max}[0,2N-3] < c^{'}=\Re\, s^{'} < 2N-1$. 
By introducing the reflection formula for the Riemann zeta function as given on p.\ 269
of Ref.\ \cite{whi73} or No.\ 9.535(3) of Ref.\ \cite{gra94}, which is
\begin{eqnarray}
2^{1-s} \,\Gamma(s) \, \zeta(s) \,\cos \Bigl( \frac{\pi s}{2}\Bigr) = \pi^s\, \zeta(1-s) \;\;,
\label{twentynine}\end{eqnarray}
we find that
\begin{align}
\sum_{k=N}^{\infty} \frac{B_{2k}}{2k(2k-1)z^{2k-1}} & = - \frac{1}{2\pi i} \int_{c-i\infty \atop{{\rm Max}[0,2N-3]<c= {\rm Re}\,s<2N-1}}^{c+i\infty}  
\;  \frac{z^{-s} \, \zeta(-s)}{\sin(s\pi)} \; ds
\nonumber\\
& + \;\; \frac{1}{2\pi i} \int_{C} \frac{1}{(2\pi z)^{2s-1}} \; \frac{\Gamma(2s-1) \zeta(2s)}{\sin(s \pi)} \; ds \;\;,
\label{thirty}\end{align}
where the prime superscripts have been dropped. For $N \!=\! 1$, Eq.\ (\ref{thirty}) reduces to
\begin{align}
\sum_{k=1}^{\infty} \frac{B_{2k}}{2k(2k-1)z^{2k-1}} & = - \frac{1}{2\pi i} \int_{c-i\infty \atop{{\rm Max}[0<c= {\rm Re}\,s<1}}^{c+i\infty}  \; 
\frac{z^{-s} \, \zeta(-s)}{\sin(s\pi)} \; ds
\nonumber\\
& + \;\; \frac{1}{2\pi i} \int_{C} \frac{1}{(2\pi z)^{2s-1}} \; \frac{\Gamma(2s-1) \zeta(2s)}{\sin(s \pi)} \; ds \;\;.
\label{thirtyone}\end{align}
This is identical to the result with which Paris begins in Ref.\ \cite{par14} except it includes the contour integral along the 
great arc. The problem is that the contour along the great arc does not always vanish and is in fact divergent for certain values of $z$. When this 
occurs, it must be reflected in the series on the lhs becoming divergent. So let us examine the series on the lhs.

By introducing Eq.\ (\ref{twentyseven}) into the lhs of Eq.\ (\ref{twentyeight}) and replacing the Riemann zeta function by its Dirichlet series
form, we obtain
\begin{eqnarray}
L(N,z)= \sum_{k=N}^{\infty} \frac{B_{2k}}{2k(2k-1)z^{2k-1}}= z \sum_{n=1}^{\infty} \sum_{k=N}^{\infty} \frac{(-1)^k \,\Gamma(2k-1)}{(2\pi n)^{2k}}\;\;.  
\label{thirtytwo}\end{eqnarray} 
The inner series looks very familiar indeed. Except for the factor of 2 next to $k$ inside the gamma function, it is basically a Type I terminant, which 
we have studied earlier. Because of this factor such terminants are called generalized terminants in Ref.\ \cite{kow09}. Specifically, a  
generalized Type I terminant is defined as
\begin{eqnarray}
S^{I}_{p,q}\Bigl( N,z^{\beta} \Bigr)= \sum_{k=N}^{\infty} \Gamma(pk+q) (-z^{\beta})^k \;\;.
\label{thirtythree}\end{eqnarray}
Therefore, the inner series in Eq.\ (\ref{thirtytwo}) can be represented by $S^{I}_{2,-1}(N,-(1/2\pi n)^2)$. Let us now introduce the integral representation
for the gamma function into Eq.\ (\ref{thirtythree}) and interchange the order of the summation and integration, which represent the first two
steps in Borel summation. As displayed on p.\ 156 of Ref.\ \cite{kow09}, we find that
\begin{eqnarray}
S^{I}_{p,q}\Bigl( N,z^{\beta} \Bigr)= \int_0^{\infty} e^{-t} \,t^{q-1} 
\sum_{k=N}^{\infty}\Bigl(-z^{\beta}t^p\Bigr)^k \;\;.
\label{thirtyfour}\end{eqnarray}
For the above integral to be convergent always, we require that $\Re\, (pN+q)>0$. Once again, we encounter the geometric series, which we have
already studied in detail. In fact, from Equivalence\ (\ref{fourteen}) we know that
\begin{eqnarray}
\sum_{k=N}^{\infty}\Bigl(-z^{\beta}t^p\Bigr)^k 
\begin{cases}
\equiv \left( -z^{\beta} t^p \right)^{N}/(1+ z^{\beta} t^p)\;\;, & \;\; \Re \, \left( -z^{\beta} t^p \right) \geq 1\;\;, \\
= \left( -z^{\beta} t^p\right)^{N}/(1+z^{\beta} t^p)\;\;, & \;\; \Re \, \left( -z^{\beta} t^{p}\right) <1.
\end{cases}
\label{thirtyfive}\end{eqnarray}
For the asymptotic series in the complete asymptotic expansion of $\ln \Gamma(z)$, viz.\ $L(N,z)$, we have see that $\beta=2$ and $p=2$. Therefore, 
the series possesses an infinity in its remainder whenever $\Re\, (-z^2t^2) >1$. As $t$ ranges from zero to infinity, this means that the integral 
will yield an infinite value for $\Re \, z^2 < 0$. Therefore, there is no doubt that the integral along the great arc must also be infinite for these 
values of $z$ or else Eq.\ (\ref{thirty}) is nonsense. That is, MB regularization amounts to the removal of the contour along the great arc, which is 
either zero or infinity. Consequently, Eq.\ (\ref{thirtyone}) can be expressed as
\begin{align}
\sum_{k=1}^{\infty} \frac{B_{2k}}{2k(2k-1)z^{2k-1}} \equiv - \frac{1}{2\pi i} \int_{c-i\infty \atop{{\rm Max}[0<c= {\rm Re}\,s<1}}^{c+i\infty}  \; 
\frac{z^{-s} \, \zeta(-s)}{\sin(s\pi)} \; ds \;\;.
\label{thirtyfivea}\end{align}
Thus, we have seen that the result beginning the analysis in Ref.\ \cite{par14} has already been regularized. More importantly, it is simply 
incorrect to say that there is no need for regularization to appear in an asymptotic series. All asymptotic series possess
regions/sectors where they are either divergent or conditionally convergent because the method used to determine them is improper as explained in Ref.\
\cite{kow09}. Regularization is necessary for correcting the divergence resulting from the application of an asymptotic method. On the other hand, we 
have seen that an asymptotic series is not divergent everywhere. It has regions where it is convergent, although such series are conditionally convergent 
rather than possessing a region of absolute convergence like the geometric and binomial series. Furthermore, the above statement is not complete because 
the values of ${\rm arg} \,z $ over which it is valid have not been stated. This leads to the concept of an asymptotic form, which is discussed in the 
following section.

\section{The Stokes Phenomenon and Smoothing}
There seems to be a line of thinking amongst the proponents of the smoothing of the Stokes phenomenon that the interpretation of the Stokes phenomenon 
given in Refs.\ \cite{kow14} and \cite{kow09} is not the same as theirs and that as a consequence, although the calculations in these references are 
correct as asserted by Paris in Ref.\ \cite{par14}, we are not talking about the same phenomenon. This section aims to refute such a proposition.

Before discussing the specific comments made by Paris concerning smoothing, let us consider the following paragraphs from Berry's famous paper on the 
asymptotic smoothing of Stokes's discontinuities \cite{ber89}. There it is written:

\indent
{\it The conventional view (Stokes 19864) is powerfully (and unconventionally) argued by Dingle in I. It asserts that the change in $S$ is discontinuous
and localized at the Stokes line: on one side, $S$ takes a value, $S_{-}$ say; on the other, $S= S_{-}+1$; on the line itself, $S=S_{-} +1/2$. For
the example (3) the intuition behind this view is illustrated by figure 1, which shows how the steepest-descent contours of the integral
(${\rm Im}\, t^2= {\rm Im}\, Z^2$) change discontinuously across the Stokes line, suddenly bringing in the subdominant contribution from the
stationary point at $t=0$ ($S_{-}=0$ in this case).  It is worth repeating .....

From the context it is clear that Stokes is referring to asymptotic series interpreted by truncation near their least term. My aim here is to dispel
Stokes's mist and show that his discontinuity is an artefact of poor resolution: with the appropriate magnification, $S$ changes smoothly.}
\noindent

From the above paragraphs it is obvious that we are not just talking about a minor difference in interpretation, but a completely different view or
understanding of the Stokes phenomenon with major ramifications to asymptotics. Moreover, Berry is aiming to demonstrate that since Stokes and Dingle 
are the main proponents of the conventional view, both of them are wrong. In fact, Paris refers to the dramatic change in the field of asymptotics as a result
of Berry's paper in Ch.\ 6 of his book \cite{par01}. Strangely, it is today referred to as Stokes smoothing when it should be the Berry smoothing of
the Stokes phenomenon since Stokes \cite{sto04} never entertained this view. As we shall see, the view is very disturbing because if it is 
true, then it would be impossible to obtain exact values of a function from its complete asymptotic expansion as has been accomplished in Refs.\ 
\cite{kow14}, \cite{kow09}, \cite{kow11b}-\cite{kow95} and \cite{kow002}.

For the benefit of the reader, the smoothing effect takes place near what are known as Stokes lines. These are lines where an asymptotic series
attains maximal dominance over all other terms in the asymptotic expansion, which appears as Rule B in Ch.\ 1 of Dingle's book \cite{din73}. Basically,
this means that Stokes lines or rays are determined by those phases in which the terms in an asymptotic series are homogeneous in phase and all of
the same sign (Rule A). The regions or sectors bordered by Stokes lines are referred to as Stokes sectors. In terms of the Type I terminant given in 
``Eqs."\ (\ref{seven}) and (\ref{eight}), the Stokes lines occur at ${\rm arg} \,x =(2k+1) \pi$, where $k$ is an integer because all the terms in the 
asymptotic series are all positive real numbers with no complex or imaginary parts at all. Since they occur virtually outside the principal branch, Dingle 
does not study the Stokes phenomenon in relation to Type I terminants in Ch.\ 22 of Ref.\ \cite{din73}. Instead, he applies the conventional view to 
Type II terminants, which are identical to Type I terminants except the factor of $(-1)^k$ is missing. For these series the Stokes lines occur at 
${\rm arg}\, x=2k \pi$ (including zero), well within the principal branch of the complex plane. In actual fact, this is a shortcoming in his great book 
because often an asymptotic series is in powers of $x^{\beta}$, where $\beta >1$. Then the conditions for the Stokes lines become either ${\rm arg}\, z= 
(2k+1) \pi/\beta$ for a generalized Type I terminant or ${\rm arg}\,z=2k \pi/\beta$ for a generalized Type II terminant. Hence there can be many Stokes 
lines lying in the principal branch, which is usually taken to be $(-\pi,\pi]$. Because of this, all phases and sectors are studied for both types of 
generalized terminants in Ref.\ \cite{kow09}.    

Not long after Berry's paper appeared, Olver \cite{olv90}  claimed that he had ``rigorously" proven that the smoothing did occur near a Stokes line and 
that Stokes had got it wrong after all. Basically, Olver re-expressed Berry's form for the Stokes multiplier of the subdominant asymptotic solution
of the one-dimensional Helmholtz equation, viz.\
\begin{eqnarray}
S_n(F) = \frac{1}{2} -\frac{i}{2 \pi} \; P \! \int_0^{\infty} dt\; \frac{t^{n-\beta}}{1-t}\, e^{F(1-t)}\;\;,
\label{thirtyfiveb}\end{eqnarray}
as
\begin{eqnarray}
S_n=-\frac{1}{2\pi i} \int_{-1}^{\infty} \frac{d\tau}{\tau} \, e^{-A(\tau-\ln(1+\tau))}\,(1+\tau)^{\mu}
\; e^{-i B\tau} \;\;,
\label{thirtyfivec}\end{eqnarray}
where $2z \!=\! -A \!-\! i B$, $n \!-\! 1 \!=\! A \!+\! \mu$, $A$ is large, real and positive, $n$ is an integer and $B$ and $\mu$ are real. 
From Laplace's method, which is a variant of the method of expanding most of the exponential mentioned earlier, one expects as $A \!\to \! \infty$ 
that the major contribution to the above integral comes from the neighbourhood around $\tau \!=\!0$. Hence Olver expands the first exponential 
in the above integral in powers of $\tau$ and extends the lower limit of the integral to $-\infty$. In so doing, he eventually obtains
\begin{align}
S_n & \overset{\sim}{=} \frac{1}{2} + \frac{1}{2} \; {\rm erf}\Bigl( \frac{B}{\sqrt{2A}}\Bigr)
+i \; \frac{e^{-B^2/2A}}{\sqrt{2 \pi A}} \left( \mu +\frac{1}{3}-\frac{B^2}{3A}+
\Bigl( \frac{1}{A} -\frac{B^2}{A^2} \Bigr) \right.
\nonumber\\
& \left. \times \;\;  \Bigl( \frac{\mu}{3} -\frac{\mu^2}{2} + \frac{\mu^3}{6} \Bigr) + \cdots \right) 
-\frac{B}{2A}\; \frac{e^{-B^2/2A}}{\sqrt{2 \pi A}} \Bigl( \mu-\mu^2+ \cdots \Bigr) \;\;. 
\label{thirtyfived}\end{align}
Hence we see the emergence of the error function appearing in the Stokes multiplier. It represents the leading order term of complicated asymptotic
expansion, which at best is conditionally convergent, but whose terms diverge in magnitude eventually. This is what is referred to as rigorous 
mathematics. All that has happened here is that more than the first few orders have been obtained via standard Poincar$\acute{\rm e}$ asymptotics with
the remaining terms denoted by + $\cdots$. On the other hand, not one asymptotic result has been truncated in Ref.\ \cite{kow14}. The only instance 
of truncation there was to a convergent series in order to expedite the calculations and even then, it was after $10^5$ terms in the series. For a more 
detailed description of the terms missing in the above result, the reader should consult Ch.\ 6 and the appendix of Ref.\ \cite{kow09}. So, it 
appears that Berry's smoothing of the Stokes phenomenon is an ``artefact of poor resolution" rather than the other way around with Stokes mist.

What Berry, Olver and Paris have missed is that Dingle came up with an explanation of the Stokes phenomenon when discussing
Borel summation of Type II terminants in Ch.\ 22 of his book. There he states that
\begin{eqnarray}
\sum_{k=n}^{\infty} \frac{\Gamma(k+\alpha +1)}{x^k}\equiv \frac{\Gamma(n+\alpha+1)}{x^n} \; {\bar \Lambda}_{n+\alpha} (-x) \;\;, \quad {\rm arg}\, x=0, 
\label{thirtyfivee}\end{eqnarray}
where 
\begin{eqnarray}
{\bar \Lambda}_s(-x) =\frac{1}{\Gamma(s+1)} \; P \int_0^{\infty} \frac{t^s \,e^{-t}}{1-t/x} \;dt \;\;,
\label{thirtysix}\end{eqnarray}
and $P$ denotes that the Cauchy principal value must be taken. In the above we have replaced the equals sign in Dingle's ``Eq."\ (16) by an equivalence 
symbol in accordance with the discussion on regularization in the previous section. The above result is actually an asymptotic form because it 
indicates the values of $x$ for which it is valid. Moreover, the result applies only to a Stokes line since according to the conventional view 
there is a jump discontinuity as $x$ moves either above or below the positive real axis. So, what is this jump discontinuity? Well, if we move 
to the next section in Dingle's book, then we see the jump discontinuities emerge due to semi-circular residues of the integral in Eq.\ 
(\ref{thirtysix}). Therefore, we are basically interpreting the integral as a Cauchy integral with a line contour along the positive real axis, but as 
$x$ moves either above or below the positive axis, we  have to take into account the semi-residues of the Cauchy integral. Dingle states for 
${\rm arg}\,x>0$, the semi-residue is evaluated in an anticlockwise direction, while for ${\rm arg}\,x <0$, it is taken in a clockwise direction. 
As a consequence, he finds that
\begin{eqnarray}
{\bar \Lambda}_s(-x) =  \begin{cases} \Lambda_s(-x)  + i \pi x^{s+1} e^{-x} /\Gamma(s+1) \;, \quad 0< {\rm arg}\, x < 2 \pi \;, \cr
\Lambda_s(-x)  - i \pi x^{s+1} e^{-x} /\Gamma(s+1) \;, \quad -2\pi < {\rm arg}\, x < 0 \;.
\end{cases} 
\label{thirtyseven}\end{eqnarray}
In the above results $\Lambda_s(-x)$ is given by Eq.\ (\ref{sevena}). It should also be noted that these results are virtually identical to the 
Plemelj relations as discussed on p.\ 414 of Ref.\ \cite{car66}. 
 
We can summarize Dingle's results by expressing them as
\begin{align}
\sum_{k=n}^{\infty} \frac{\Gamma(k+\alpha +1)}{x^k}\equiv \frac{\Gamma(n+\alpha+1)}{x^n} \; \Lambda_{n+\alpha} (-x) 
+ \frac{2 \pi i\,x^{s+1} \,e^{-x}}{ \Gamma(n+\alpha+1)}\; S \;\;, 
\label{thirtyeight}\end{align} 
where we take the Cauchy principal value when evaluating $\Lambda_{n+\alpha}(-x)$ for ${\rm arg}\, x=0$ and the factor $S$ 
multiplying the full residue is given by
\begin{eqnarray} 
S= \begin{cases}  1/2\;\;, \quad 0< {\rm arg}\,x <2\pi \;\;, \cr 
0\;\;, \quad {\rm arg}\, x =0 \;\;,\cr
-1/2 \;\;, \quad -2\pi < {\rm arg}\, x<0 \;\;.
\end{cases}
\label{thirtynine}\end{eqnarray}
This indeed looks familiar. It is in fact the conventional view of the Stokes multiplier as described by Berry at the beginning of this 
section with $S_{-} \!=\! -1/2$. Thus, Dingle has provided us with the key to the Stokes phenomenon. We simply interpret the integrals
obtained by Borel summation as Cauchy integrals and then investigate their singular behaviour. Moreover, from these results we see
that it is necessary that the values of ${\rm arg}\, z$ need to be specified for an asymptotic series. The combination of the
asymptotic expansion with the values of ${\rm arg}\, z$ over which it is valid gives rise to asymptotic forms. 

As indicated earlier, the asymptotic series in the Stirling's approximation for $\ln \Gamma(z)$ involves Type I terminants, which are 
not discussed in Ch.\ 22 of Dingle's book \cite{din73}. For this situation we require the material in Ch.\ 10.1 of Ref.\ \cite{kow09},
which deals with the Borel summation of generalized Type I terminants given by Eq.\ (\ref{thirtythree}). In fact, Eq.\ (\ref{thirtythree}) is
extended by introducing a factor of $\exp(-2li \pi)$ with $z^{\beta}$. Thus, Borel summation of a generalized Type I terminant yields    
\begin{eqnarray}
S^{I}_{p,q} \! \left( N,z^{\beta} e^{-2 li \pi} \right) \equiv  (-1)^N p^{-1}  z^{\beta (N-1)} \int_C 
\frac{s^{N+q/p-1}\, e^{-s^{1/p}}}{s-(-z^{-\beta} e^{2 l i \pi})} \; ds \;\; .
\label{forty}\end{eqnarray}
If we let 
\begin{eqnarray}
f(s) = z^{\beta(N-1)} s^{N+q/p-1} \exp(-s^{1/p})/p\;\;,
\label{fortya}\end{eqnarray} 
then the rhs in the above result can be regarded as a Cauchy integral whose contour $C$ is the line contour along the positive real axis.
Furthermore, it possesses a singularity at $s \!=\!  -z^{-\beta}$. Although $\exp(2l i \pi)$ is equal to unity for all values of $l$, it was found 
in Ch.\ 7 of Ref.\ \cite{kow09} to have an effect on the MB-regularized value of $S^{I}_{p,q}(N,z^{\beta}e^{-2li \pi})$. In particular, the 
difference between the MB-regularized values of $S^I_{p,q} \!\left( N,z^{\beta} e^{-2 l i \pi}\right)$ and 
$S^I_{p,q} \!\left( N,z^{\beta}e^{-2 (l \!-\! 1) i \pi}\right)$ was found to be given by
\begin{align}
& \Delta S^{I}_l(N,z^{\beta}) = S^{I}_{p,q}  \bigl( N,z^{\beta} e^{-2 li \pi} \bigr)- S^{I}_{p,q}  
\bigl( N,z^{\beta} e^{-2 (l-1) i \pi} \bigr)  
\nonumber\\
& \equiv  \int_{c-i \infty \atop{{\rm Max}[N-1,-q/p] <c=\Re\,s<N}}^{c+i \infty}  \!\!\!\! 
 z^{\beta s} \; \Gamma(ps+q) \, e^{-(2l-1) i \pi s} \; ds  \;\; ,
\label{fortyone}\end{align}
where the lower bound on the offset $c$ is adjusted to ensure that the poles for the MB integral remain to the right of $-q/p$, if it should 
be greater than $N \!-\! 1$, in accordance with the conditions on $f(k)$ in Eq.\ (\ref{sixteena}). By making the change of variable, 
$y \!=\! ps \!+\! q$, we find that Equivalence\ (\ref{fortyone}) becomes
\begin{align}
& \Delta S^{I}_l(N,z^{\beta}) \equiv p^{-1} (z^{-\beta} \, e^{(2l-1)i \pi})^{q/p} \int_{c-i \infty \atop{c=\Re\,s>-q/p}}^{c+i \infty}  
\!\!\!\! 
\left( z^{-\beta} e^{(2l-1)i \pi} \right)^{\! -y/p} \Gamma(y)\; dy  \;\;.
\label{fortytwo}\end{align} 
The MB integral can be regarded as the inverse Mellin transform of $\exp(-x)$ with $x \!=\! z^{-\beta/p}\exp((2l-1)i \pi/p)$ in Eq.\ 
(\ref{twenty}). Consequently, the above equivalence reduces to
\begin{eqnarray}
\Delta S^{I}_l(N,z^{\beta}) \equiv \frac{2 \pi i}{p} \;z^{-\beta q/p}\, e^{(2l-1)i q\pi/p} 
\exp \!\left(-z^{-\beta/p} e^{(2l-1)i \pi/p} \right) \;\; . 
\label{fortythree}\end{eqnarray}

On p.\ 412 of Ref.\ \cite{car66} we see that the Cauchy integral on the rhs of Equivalence\ (\ref{forty}) develops jump discontinuities as 
$-z^{-\beta} e^{2i l \pi }$ moves across the line contour. This means that while $z^{-\beta} e^{2li \pi }$ is located in a branch of the complex 
plane, say $(2j \!-\!1)\pi \!<\! {\rm arg}\, (z^{-\beta}e^{2l i  \pi }) \!<\! (2j \!+\! 1) \pi$, the regularized value is given by the Cauchy 
integral, but once $z^{-\beta}e^{2li  \pi }$ moves outside of this branch, it acquires extra terms or else the regularized value would be the 
same for all Stokes sectors. Hence the regularized value of $S^{I}_{p,q}(N,z^{\beta}e^{-2 li \pi })$ cannot be represented solely by a 
Cauchy integral. As a result, another problem emerges. We need to determine the specific Stokes sector over which the Cauchy 
integral is the sole contribution to the regularized value. This appears to be arbitrary, much like selecting a principal branch in the complex 
plane. Thus, a primary Stokes sector must be nominated. This is taken to be the $j \!=\!l \!=\! 0$ branch of the complex plane. Consequently, the 
Cauchy integral on the rhs of Equivalence\ (\ref{forty}) represents the regularized value only for $-\pi/\beta \!<\! {\rm arg}\,z \!<\! \pi/\beta$. 
That is, the $j \!=\! l \!=\! 0$ branch reduces to the principal branch of the complex plane when $\beta \!=\! 1$, thereby yielding the regularized 
value of the first type of terminant as given on p.\ 406 of Dingle's book \cite{din73}. Consequently, for $l \!=\!0$ we arrive at
\begin{eqnarray}
S^{I}_{p,q} \! \left( N,z^{\beta} \right) \equiv  (-1)^N p^{-1}  z^{\beta (N-1)} \int_C  
\frac{s^{N+q/p-1}\, e^{-s^{1/p}}}{s-(-z^{-\beta})} \; ds \;\; ,
\label{fortyfour}\end{eqnarray} 
which is only valid for $-\pi/\beta \!<\! {\rm arg}\, z \!<\! \pi/\beta$.

We are now in a position to determine the jump discontinuities that apply over the secondary Stokes sectors in the complex plane. 
First, we note that the Cauchy integral in Equivalence\ (\ref{forty}) is singular whenever ${\rm arg}\,z^{-\beta} \!=\! (2 j \!+\! 1)\pi$, 
and $j$ is an integer. As $z^{-\beta} \exp((2l \!-\!1)i\pi)$ crosses from the primary sector Stokes sector to the adjacent sector or 
$j \!=\! l \!+\! 1$ branch of the complex plane, $-z^{-\beta}e^{2 li \pi}$ moves from below the line contour or positive real axis to above 
the axis. During this transition the Cauchy integral becomes undefined when $-z^{-\beta}e^{2li \pi}$ is situated on the positive real axis. 
To evaluate the residue, let us consider an infinitesimal circle around the pole at $s \!=\! -z^{-\beta}\exp(-2 li \pi )$. 
At this stage, we shall not be concerned with whether we are considering a complete rotation or a semi-circular rotation around the pole. 
Nor will we be concerned with the direction of the indentation. Hence we shall assume that the infinitesimal indentation begins 
at an angle, $\gamma_1$ in the complex plane, and ends at another angle, $\gamma_2$. Then we find that the contribution from the 
pole from the Cauchy integral in Equivalence\ (\ref{forty}) is given by  
\begin{align}
& I^{I} =  i p^{-1} (-1)^N z^{\beta(N-1)} \lim_{\epsilon \to 0} \int_{\gamma_1}^{\gamma_2}  
\left(z^{-\beta} e^{i(2l-1) \pi} \right)^{\! N+q/p-1}   
\nonumber \\ 
&  \left. \times \;\; \exp \!\left(-z^{-\beta/p}e^{i(2l-1)\pi/p} \right)\; d\gamma
= i (-1)^N \Delta \gamma \, f(s) \right|_{s=z^{-\beta}\exp((2l-1)i\pi)} \;\;,
\label{fortyfive}\end{align}
where $\Delta \gamma \!=\! \gamma_2 \!-\!\gamma_1$. Therefore, we have seen that $f(s)$ or Eq.\ (\ref{fortya}) with $s \!=\! -z^{-\beta} \exp((2l-1)i \pi)$ 
represents the residue of the Cauchy integral. Moreover, if $\Delta \gamma \!=\! 2 \pi$, which corresponds to a complete rotation in an 
anticlockwise direction, then we see that the residue of the Cauchy integral given by Equivalence\ (\ref{forty}) is identical to the difference 
of the regularized values for $S^I_{p,q}(N,z^{\beta}e^{-2li \pi})$ and $S^I_{p,q}(N,z^{\beta}e^{-2(l-1)i\pi})$ or Equivalence\
(\ref{fortythree}). Therefore, the above result confirms the remarkable insight made by Dingle on p.\ 412 of Ref.\ \cite{din73} that the 
jump discontinuity due to crossing Stokes sectors is dependent upon the singular behaviour of the Cauchy integral that emerges from the 
introduction of the regularized value of the geometric series during Borel summation.

By putting $l \!=\! 1$ in Equivalence\ (\ref{fortythree}), we arrive at 
\begin{align}
S^{I}_{p,q} \! & \left( N,z^{\beta}e^{-2i \pi} \right)-S^{I}_{p,q} \! \left( N,z^{\beta} \right)
\equiv \frac{2\pi i}{p} \; z^{-\beta q/p} 
e^{iq \pi/p} \exp \!\left( -z^{-\beta/p}e^{i \pi/p} \right) \;\; .
\label{fortysix}\end{align}
Since the above equivalence is valid only for $-\pi/\beta \!<\! {\rm arg}\, z \!<\! \pi/\beta$, we can replace $S^{I}_{p,q}(N,z^{\beta})$ 
by its regularized value as given by Equivalence\ (\ref{fortyfour}). Hence Equivalence\ (\ref{fortysix}) becomes
\begin{align}
S^{I}_{p,q} \!  \left( N,  z^{\beta}e^{-2i \pi} \right) & \equiv  (-1)^N  p^{-1}\, z^{\beta (N-1)} 
\int_C  \frac{s^{N+q/p-1}\; e^{-s^{1/p}}}{s-(-z^{-\beta})} \; ds  
\nonumber\\
& + \;\; \frac{2\pi i}{p} \; z^{-\beta q/p}  e^{iq \pi/p} \exp \!\left( -z^{-\beta/p}e^{i \pi/p} 
\right) \;\; .
\label{fortyseven}\end{align}
Now we replace $ z\exp(-2i\pi/\beta)$ by $z_{*}$, which yields
\begin{align}
S^{I}_{p,q} \!  \left( N, z_{*}^{\beta} \right) & \equiv  (-1)^N  p^{-1}\, z^{\beta (N-1)} 
\int_C  \frac{s^{N+q/p-1}\; e^{-s^{1/p}}}{s+z^{-\beta}} \; ds 
\nonumber\\
& + \;\; \frac{2\pi i}{p} \; z^{-\beta q/p}  e^{iq \pi/p} \exp \!\left( -z^{-\beta/p}e^{i \pi/p} 
\right) \;\; ,
\label{fortyeight}\end{align}
where $-3 \pi/\beta \!<\! {\rm arg}\, z_{*} \!<\! -\pi/\beta$ and $-\pi/\beta \!<\! {\rm arg}\,z \!<\! 
\pi/\beta$. The terms on the rhs of Equivalence\ (\ref{fortyeight}) have been left in terms of 
$z$ to emphasize the fact that when they are evaluated, they are done so in the primary Stokes sector. 
That is, the regularized value of the series on the lhs of Equivalence\ (\ref{fortyeight}), which applies for
values of $z_{*}$ lying in the sector of $(-3\pi/\beta,-\pi/\beta)$, is determined by evaluating the terms
on the rhs for the corresponding values of $z_{*}$ lying in the primary Stokes sector. This anomaly arises 
from the fact that if software packages such as Mathematica are used to carry out calculations of the 
rhs of the above equivalence in determining regularized value of the generalized Type I terminant outside 
the primary sector, then they will only do so for values of the complex variable lying in the principal branch 
of the complex plane. That is, in order to obtain regularized value outside the primary sector, we need 
to evaluate forms where the complex variable lies inside it. 

To make Equivalence\ (\ref{fortyeight}) appear less awkward, we replace $z_{*}$ and $z$, respectively 
by $z$ and $z_1$, where the latter is defined as $z_1 \!=\! z \exp(2 i \pi/\beta)$. Then the above 
equivalence can be re-written as
\begin{align}
S^{I}_{p,q} \!  \left( N, z^{\beta} \right) & \equiv  (-1)^N  p^{-1}\, z_1^{\beta (N-1)} 
\int_C  \frac{s^{N+q/p-1}\; e^{s^{1/p}}}{s-(-z_1^{-\beta})} \; ds 
\nonumber\\
& + \;\; \frac{2\pi i}{p} \; z_1^{-\beta q/p}  e^{iq \pi/p} \exp \!\left( -z_1^{-\beta/p}e^{i \pi/p} 
\right) \;\; ,
\label{fortyeighta}\end{align}
where $-3 \pi/\beta \!<\! {\rm arg}\, z \!<\! -\pi/\beta$.  It is this result that is used to derive the 
regularized value of $S^{I}_{p,q}(N,z^{\beta})$ for any Stokes sector as a result of the continuous 
rotation of $z^{\beta}$ in the complex plane. 

So far, we have only been concerned with Stokes sectors, but the behaviour on the Stokes lines, where the singularity 
of the Cauchy integral is situated, is also important for developing a complete understanding of the Stokes phenomenon. 
To obtain the Borel-summed regularized value of $S^{I}_{p,q} (N,z^{\beta} )$ on the Stokes line given by ${\rm arg}\,z \!=\! -\pi/\beta$, 
we need to invoke Rule 8a presented in Ch.\ 3 of Ref.\ \cite{kow09}. The various rules appearing there are virtually a re-expression 
of those appearing in Ch.\ 1 of Dingle's book \cite{din73}. Rule 8a states that the asymptotic form along a Stokes line is 
the average of the asymptotic forms in the adjoining Stokes sectors except that the Cauchy principal value must be evaluated. 
Although there is no statement about the Cauchy principal value in Dingle's rules, we have seen that it is implied as 
demonstrated by ``Eq."\ (\ref{thirtysix}). Therefore, the Borel-summed regularized value of a generalized Type I terminant for
${\rm arg}\, z \!=\! -\pi/\beta$ is given by
\begin{align}
S^{I}_{p,q} \!  \left( N, z \right) & \equiv  -  p^{-1}\, |z|^{\beta (N-1)} 
P \int_0^{\infty} \frac{t^{N+q/p-1}\; e^{-t^{1/p}}}{t-|z|^{-\beta}} \; dt  
\nonumber\\
& + \;\; \frac{\pi i}{p} \; |z|^{-\beta q/p}  \exp \!\left( -|z|^{-\beta/p} \right) \;\; .
\label{fortynine}\end{align}
From this result we observe that the regularized value is composed of the regularized value in the primary 
Stokes sector except the principal value of the Cauchy integral is now evaluated and only the 
semi-residue contribution is taken in a anticlockwise direction ($\Delta \gamma \!=\! \pi$) of the $l \!=\! 1$ 
version of Eq.\ (\ref{fortyeighta}).
 
If we put $l \!=\!2$ in Equivalence\ (\ref{fortyone}), then with the aid of Equivalence\ (\ref{fortythree}) we find 
that
\begin{align}
S^{I}_{p,q} \!\left(N,z^{\beta}e^{-4i \pi} \right) & -S^{I}_{p,q} \!\left(N,z^{\beta} e^{-2i \pi} 
\right) \equiv \frac{2\pi i}{p} \; z^{-\beta q/p} 
e^{3iq \pi/p} \exp \!\left(-z^{-\beta/p} e^{3i q\pi/p} \right) \;\;.
\label{fifty}\end{align}
We can replace $S^{I}_{p,q}(N,z^{\beta}\exp(-2i \pi))$ by introducing Equivalence\ (\ref{fortyseven}) 
into the above result. This yields
\begin{align}
& S^{I}_{p,q} \!\left(N,z^{\beta}e^{-4i \pi} \right) \equiv (-1)^N  p^{-1}\, z^{\beta (N-1)} 
\int_C  \frac{s^{N+q/p-1}\; e^{-s^{1/p}}}{s-(-z^{-\beta})} \;ds  
+ \frac{2\pi i}{p} \; z^{-\beta q/p}  
\nonumber\\
& \times  e^{3iq \pi/p} \exp \!\left(-z^{-\beta/p} e^{3i q\pi/p} \right) +  \frac{2\pi i}{p} 
\; z^{-\beta q/p}  e^{iq \pi/p} \exp \!\left( -z^{-\beta/p}e^{i \pi/p} \right)  .
\label{fiftyone}\end{align}
Next we replace $z\exp(-4i \pi/\beta)$ by $z$ on the lhs, while on the rhs $z$ is replaced by $z_2$,
the latter now being equal to $z\exp(4i \pi/\beta)$. Consequently, we arrive at 
\begin{align}
S^{I}_{p,q} \!\left(N,z^{\beta} \right)  & \; \equiv (-1)^N  p^{-1}\, z_{2}^{\beta (N-1)} 
\int_C  \frac{s^{N+q/p-1}\; e^{-s^{1/p}} }{s+z_{2}^{-\beta}} \; ds  
+ \frac{2\pi i}{p} \; z_2^{-\beta q/p} 
\nonumber\\
& \times \sum_{j=1}^2 e^{i(2j-1)q \pi/p} \exp \!\left(-z_2^{-\beta/p} e^{i(2j-1) \pi/p} \right) \;\;  ,
\label{fiftytwo}\end{align}
where $-(2(2) \!+\! 1)\pi/\beta \!<\! {\rm arg}\, z \!<\!-(2(2) \!-\!1)\pi/\beta$ or $-5\pi/\beta \!<\!
{\rm arg}\, z \!<\! -3\pi/\beta$. When $z$ lies on the Stokes line that borders the $l \!=\!1$ and 
$l \!=\! 2$ sectors, i.e.\ where ${\rm arg}\,z \! =\! -3\pi/\beta$, all we need to do is average Equivalences\ 
(\ref{fiftyone}) and (\ref{fiftytwo}) and take the Cauchy principal value of the resulting 
contour integral. Before we can average the two equivalences we must replace $z_1$ in Equivalence\
(\ref{fiftyone}) by $z \exp(2i \pi/\beta)$ and $z_2$ in Equivalence\ (\ref{fiftytwo}) by 
$z\exp(4i \pi/\beta)$ so that $z$ is the same variable in both equivalences. Then taking the average
of the two modified equivalences and setting ${\rm arg}\, z$ equal to $-3 \pi/\beta$, we find that the 
regularized value of the generalized Type I terminant on the secondary Stokes line is given by
\begin{align}
 S^{I}_{p,q} \! \left( N, z^{\beta} \right)& \;   \equiv  -  p^{-1}\, |z|^{\beta (N-1)} 
P \int_0^{\infty} \frac{t^{N+q/p-1}\;e^{-t^{1/p}}}{t-|z|^{-\beta}} \; dt   
+ \frac{2\pi i}{p} \; |z|^{-\beta q/p} e^{2iq\pi/p} 
\nonumber\\
& \times \;\; \exp \!\left( -|z|^{-\beta/p}e^{2i \pi/p} \right)  +  \frac{\pi i}{p} \, 
|z|^{-\beta q/p}  \exp \!\left( -|z|^{-\beta/p} \right) \;\;   .
\label{fiftythree}\end{align}

A pattern is now emerging that enables us to determine the regularized value of a generalized Type I 
terminant for the Stokes sectors due to clockwise rotations of $z$. We simply replace the upper 
limit 2 by $M$ in the sum on the rhs of Equivalence\ (\ref{fiftytwo}). As a consequence, we find that
the regularized value of $S^{I}_{p,q}(N,z^{\beta})$ for  $-(2M \!+\! 1)\pi/\beta \!<\! 
{\rm arg}\, z \!<\! -(2M \!-\! 1)\pi/\beta$ is given by
\begin{align}
S^{I}_{p,q}(N,z^{\beta}) & \;\equiv (-1)^N z_{M}^{\beta(N-1)} p^{-1} \int_{0}^{\infty}  \frac{t^{N+q/p-1}\, 
e^{-t^{1/p}}}{t+z_{M}^{-\beta}} \; dt+ \frac{2 \pi i}{p}\; z_M^{-\beta q/p}
\nonumber\\
& \times  \;\; \sum_{j=1}^{M}\,e^{i(2j-1) q \pi /p}\, \exp \!\left( -z_M^{-\beta/p} e^{i(2j-1) \pi/p} 
\right) \;\;, 
\label{fiftyfour}\end{align}
where $z_{M} \!=\! z \exp(2Mi \pi/\beta)$. This equivalence, which appears as (32) in the beginning of 
the proof to Thm.\ 2.1 in Ref.\ \cite{kow14}, represents the base or platform from which the exactification of Stirling's 
approximation can be carried out for positive values of ${\rm arg}\, z$. It should also be noted that negative values of 
${\rm arg}\,z$ correspond to positive values of ${\rm arg}\, z$ in Stirling's approximation since the asymptotic series 
in the latter, viz.\ $S(z)$, is in inverse powers of $z$. For the Stokes line given by ${\rm arg}\,z \!=\!-(2M \!+\!1)
\pi/\beta$, one simply determines the average of the regularized values for the two adjacent Stokes sectors bordered 
by the line and evaluates the principal value of the Cauchy integral. After a little algebra one finds that 
the regularized value of the generalized Type I terminant on the Stokes line reduces to   
\begin{align}
&  S^{I}_{p,q}(N,z^{\beta}) \equiv -|z|^{\beta(N-1)} p^{-1} \; P \int_{0}^{\infty}  \frac{t^{N+q/p-1}\, 
e^{-t^{1/p}}} {t-|z|^{-\beta}} \; dt +\frac{2 \pi i}{p}\; |z|^{-\beta q/p}
\nonumber\\
& \times \;\;\sum_{j=1}^{M}e^{2ij q \pi /p} \exp \!\left( -|z|^{-\beta/p} e^{2ij \pi/p} \right) 
 + \frac{\pi i}{p} \; |z|^{-\beta q/p}\, \exp \!\left( -|z|^{-\beta/p} \right) \;\;.
\label{fiftyfive}\end{align}
This equivalence appears as (45) in the proof of Thm.\ 2.1 in Ref.\ \cite{kow14}.

For negative values of ${\rm arg}\, z^{-\beta}$ in the complex plane, the rotations of $z^{-\beta}$ are
clockwise. To obtain the regularized value of a generalized Type I terminant for this case, we put $l \!=\!0$ 
in Equivalences\ (\ref{fortyone})-(\ref{fortythree}), which yields
\begin{align}
S^{I}_{p,q} \! & \left( N,z^{\beta}e^{2i \pi}\right)- S^{I}_{p,q} \! \left( N,z^{\beta} \right)
\equiv - \frac{2\pi i}{p} \; z^{-\beta q/p} 
\nonumber\\
& \hspace{1 cm} \times \;\; e^{-iq \pi/p} \exp \!\left( -z^{-\beta/p}e^{-i \pi/p} \right) \;\; .
\label{fiftysix}\end{align}
Hence we see that the regularized value for the lower Stokes sector given by $S^I_{p,q}(N,z^{\beta})$
is related to the regularized value for the Stokes sector immediately above, viz.\
$S^{I}_{p,q}(N,z^{\beta}$ $ \exp(2i \pi))$, plus the residue contribution of the Cauchy integral taken in a
clockwise direction. In fact, the regularized values of all lower Stokes sectors are related to the
regularized values for the Stokes sectors immediately above them plus the residue contributions of the Cauchy
integral taken in a clockwise direction. Again, this is consistent with the conventional view of the Stokes
phenomenon. Furthermore, the above result represents the complex conjugate of Equivalence\ (\ref{fortysix}). Since 
Equivalence\ (\ref{fiftysix}) is only valid over the primary Stokes sector, i.e.\ for $-\pi/\beta \!<\! {\rm arg}\, 
z \!<\! \pi/\beta$, we can introduce the regularized value of $S^{I}_{p,q}(N,z^{\beta})$ from Equivalence\ 
(\ref{fortyfour}) into it. Then by replacing $z\exp(2i \pi/\beta)$ in the resulting equivalence with $z$, one 
obtains
\begin{align}
S^{I}_{p,q} \!  \left( N,  z \right) & \equiv  (-1)^N  p^{-1}\, z_{-1}^{\beta (N-1)} 
\int_C ds\; \frac{s^{N+q/p-1}}{s-(-z_{-1}^{-\beta})}\; e^{-s^{1/p}}  
\nonumber\\
& - \;\; \frac{2\pi i}{p} \; z_{-1}^{-\beta q/p} \, e^{-iq \pi/p} \exp \!\left( -z_{-1}^{-\beta/p}
e^{-i \pi/p} \right) \;\; ,
\label{fiftyseven}\end{align}
where $\pi/\beta \!<\! {\rm arg}\,z \!<\! 3 \pi/\beta$ and $z_{-1} \!=\! z\exp(-2i \pi/\beta)$.

For the Stokes line of ${\rm arg}\, z \!=\! \pi/\beta$, we first re-write the rhs of Equivalence\ (\ref{fiftyseven}) 
in terms of $z$. Then we average the resulting equivalence with Equivalence\ (\ref{fortyfour}). Next
we take the Cauchy principal value of the resulting integral in accordance with Rule\ 8a in Ch.\ 3
of Ref.\ \cite{kow09}. Finally, we replace $z$ by $|z| \exp(i\pi/\beta)$. Thus, we arrive at the Borel-summed 
regularized value of a generalized Type I terminant for ${\rm arg}\, z \!=\! \pi/\beta$, which is
\begin{align}
S^{I}_{p,q} \!  \left( N, z \right) & \equiv  -  p^{-1}\, |z|^{\beta (N-1)} 
P \int_0^{\infty}  \frac{t^{N+q/p-1}\;e^{-t^{1/p}}}{t-|z|^{-\beta}} \; dt 
\nonumber\\
& - \;\; \frac{\pi i}{p} \; |z|^{-\beta q/p}   \exp \!\left( -|z|^{-\beta/p} \right) \;\; .
\label{fiftyeight}\end{align}
                                                      
When $l \!=\!-1$, the combination of Equivalences\ (\ref{fortyone}) and (\ref{fortythree}) yields
\begin{align}
S^{I}_{p,q} \! & \left( N,z^{\beta}e^{4i \pi}\right)- S^{I}_{p,q} \! \left( N,z^{\beta}e^{2i \pi} \right)
\equiv - \frac{2\pi i}{p} \; z^{-\beta q/p} 
\nonumber\\
& \hspace{1 cm} \times \;\; e^{-3iq \pi/p} \exp \!\left( -z^{-\beta/p}e^{-3i \pi/p} \right) \;\; .
\label{fiftynine}\end{align}
We can express $S^I_{p,q}(N,z^{\beta}\exp(2i \pi))$ in terms of $S^{I}_{p,q}(N,z^{\beta})$ by putting
$l \!=\!0$ in Equivalences\ (\ref{fortyone}) and (\ref{fortythree}). This gives
\begin{align}
& S^{I}_{p,q} \!  \left( N,z^{\beta}e^{4i \pi} \right) \equiv S^{I}_{p,q} \! \left( N,z^{\beta} \right)
- \frac{2\pi i}{p} \; z^{-\beta q/p} 
\nonumber\\
& \times \;\; \sum_{j=1}^{2} e^{-i(2j-1)q \pi/p} \exp \!\left( -z^{-\beta/p} e^{-i(2j-1) \pi/p} \right) \;\; ,
\label{fiftyninea}\end{align}
where $(2(2) \!-\! 1) \pi/\beta \!<\! {\rm arg}\, z \!<\! (2(2)\!+\! 1)\pi/\beta$. Setting $z$ equal to
$z \exp(-4i \pi/\beta)$ and introducing the appropriate version of Equivalence\ (\ref{forty}), one eventually 
obtains
\begin{align}
S^{I}_{p,q} \!\left(N,z^{\beta} \right)  & \; \equiv (-1)^N  p^{-1}\, z_{-2}^{\beta (N-1)} 
\int_C  \frac{s^{N+q/p-1}\; e^{-s^{1/p}}}{s+z_{-2}^{-\beta}} \; ds 
- \frac{2\pi i}{p} \; z_{-2}^{-\beta q/p} 
\nonumber\\
& \times \sum_{j=1}^2 e^{-i(2j-1)q \pi/p} \exp \!\left(-z_{-2}^{-\beta/p} e^{-i(2j-1) q\pi/p} 
\right) \;\;  ,
\label{sixty}\end{align}
which, as expected, is the complex conjugate of Equivalence\ (\ref{fiftytwo}). For ${\rm arg}\,z 
\!=\! 3\pi/\beta$, we first express the rhs's of Equivalences\ (\ref{fiftyseven}) and
(\ref{sixty}) in terms of $z$ rather than $z_{-1}$ and $z_{-2}$. Then we average the resulting
equivalences and evaluate the Cauchy principal value of the resulting contour integral. Alternatively, we
can take the complex conjugate of Equivalence\ (\ref{fiftythree}). Hence we find that
\begin{align}
& S^{I}_{p,q} \! \left( N, z^{\beta} \right)  \equiv  -  p^{-1}\, |z|^{\beta (N-1)} 
P \int_0^{\infty}  \frac{t^{N+q/p-1}\; e^{-t^{1/p}}}{t-|z|^{-\beta}} \; dt 
- \frac{2\pi i}{p} \; |z|^{-\beta q/p} 
\nonumber\\
& \times \; \; e^{-2iq\pi/p} \exp \!\left( -|z|^{-\beta/p}e^{-2i \pi/p} \right)  -  \frac{\pi i}{p} \, 
|z|^{-\beta q/p}   \exp \!\left( -|z|^{-\beta/p} \right) \;\;  .
\label{sixtyone}\end{align}

From the above we see that a similar pattern is emerging as in the case of the anti-clockwise rotations
of $z^{-\beta}$. Therefore, by extending Equivalences\ (\ref{sixty}) and (\ref{sixtyone}) we are able to 
derive general forms for the Borel-summed regularized values of a generalized Type I terminant for any Stokes 
sector or line, where ${\rm arg}\,z \!>\! 0$. In particular, the generalization of Equivalence\ (\ref{sixty}) to 
$(2M \!-\! 1) \pi/\beta \!<\! {\rm arg}\, z \!<\! (2M \!+\! 1) \pi/\beta$, where $M \!>\! 0$, can be carried 
out simply by replacing 2 with $M$. Thus, we arrive at
\begin{align}
& S^{I}_{p,q}(N,z^{\beta})  \equiv (-1)^N z_{-M}^{\beta(N-1)} \, p^{-1} \int_{0} ^{\infty} \frac{t^{N+q/p-1}\, 
e^{-t^{1/p}}}
{t+z_{-M}^{-\beta}} \; dt -  \frac{2 \pi i}{p}\; z_{-M}^{-\beta q/p}
\nonumber\\
& \;\; \times \;\; \sum_{j=1}^M e^{-i(2j-1) q \pi/p}\, \exp \!\left( -z_{-M}^{-\beta/p} e^{-i(2j-1) \pi/p} 
\right) \;\;. 
\label{sixtytwo}\end{align}
In the above result $z_{-M} \!=\! z \exp(-2 Mi \pi/\beta)$. Hence we see that Equivalence\
(\ref{sixtytwo}) represents the complex conjugate of Equivalence\ (\ref{fiftyfour}). This is essentially
the equivalence given as (41) in the proof of Thm.\ 2.1 of Ref.\ \cite{kow14}. For the Stokes
line where ${\rm arg}\, z \!=\! (2M \!+\! 1) \pi/\beta$, the generalization of Equivalence\ (\ref{sixtyone})
yields
\begin{align}
S^{I}_{p,q}& (N,z^{\beta}) \equiv -|z|^{\beta(N-1)} \, p^{-1} P\int_{0}^{\infty}  \frac{t^{N+q/p-1}\, 
e^{-t^{1/p}}}{t-|z|^{-\beta}}\; dt -  \frac{2 \pi i}{p}\; |z|^{-\beta q/p} \sum_{j=1}^{M} e^{-2ij q \pi /p} 
\nonumber\\
& \;\; \times \;\;  \exp \!\left( -|z|^{-\beta/p} e^{-2ij  \pi/p} \right) -
\frac{ \pi i}{p}\; |z|^{-\beta q/p}\,  \exp \!\left( -|z|^{-\beta/p} \right) \;\; ,
\label{sixtythree}\end{align}
which represents the complex conjugate of Equivalence\ (\ref{fiftyfive}). This equivalence appears as the lower-signed
version of (45) in the proof of Thm.\ 2.1 of Ref.\ \cite{kow14}. 

At no stage in the preceding analysis has there been a requirement to introduce the concept of smoothing in the vicinity
of a Stokes line. To observe the behaviour near a Stokes line, let us consider the vicinity of the line where ${\rm arg}\,z
=-(2M+1) \pi/\beta$. For $-(2M+1) \pi/\beta <{\rm arg}\,z< -(2M-1)\pi/\beta$, Equivalence\ (\ref{fiftyfour}) yields
\begin{align}
S^{I}_{p,q}(N,z^{\beta}) & \;\equiv (-1)^N z^{\beta(N-1)} p^{-1} \int_{0}^{\infty}  \frac{t^{N+q/p-1}\, 
e^{-t^{1/p}}}{t+z^{-\beta}} \; dt+ \frac{2 \pi i}{p}\; z^{-\beta q/p}
\nonumber\\
& \times  \;\; \sum_{j=0}^{M-1}\,e^{(2j+1)i q \pi /p}\, \exp \!\left( -z^{-\beta/p} e^{(2j+1)i \pi/p} 
\right) \;\;, 
\label{sixtyfour}\end{align}
while for the next lower Stokes sector or $-(2M+3)\pi/\beta < {\rm arg}\,z <-(2M+1) \pi/\beta$, it gives  
\begin{align}
S^{I}_{p,q}(N,z^{\beta}) & \;\equiv (-1)^N z^{\beta(N-1)} p^{-1} \int_{0}^{\infty}  \frac{t^{N+q/p-1}\, 
e^{-t^{1/p}}}{t+z^{-\beta}} \; dt+ \frac{2 \pi i}{p}\; z^{-\beta q/p}
\nonumber\\
& \times  \;\; \sum_{j=0}^{M}\,e^{(2j+1)i q \pi /p}\, \exp \!\left( -z^{-\beta/p} e^{(2j+1)i \pi/p} 
\right) \;\;. 
\label{sixtyfive}\end{align}
For ${\rm arg}\,z =- (2M+1)\pi/\beta$, Equivalence\ (\ref{fiftyfive}) yields
\begin{align}
 S^{I}_{p,q}(N,z^{\beta}) & \equiv (-1)^N z^{\beta(N-1)} p^{-1} \; P \int_{0}^{\infty}  \frac{t^{N+q/p-1}\, 
e^{-t^{1/p}}} {t+z^{-\beta}} \; dt +\frac{2 \pi i}{p}\; |z|^{-\beta q/p} 
\nonumber\\
& \times \;\; \sum_{j=0}^{M-1}e^{-(2j+1)i q \pi /p} 
\exp \!\left( -z^{-\beta/p} e^{-(2j+1)i \pi/p} \right) 
 + \frac{\pi i}{p} \; z^{-\beta q/p}\, e^{-(2M+1) i q\pi/p)} \,
\nonumber\\
& \times \;\; \exp \!\left( -z^{-\beta/p} e^{-(2M+1) i \pi/p}  \right) \;\;.
\label{sixtysix}\end{align}
If we let
\begin{align}
F(z) & = (-1)^N z^{\beta(N-1)} p^{-1} \int_{0}^{\infty}  \frac{t^{N+q/p-1}\, 
e^{-t^{1/p}}}{t+z^{-\beta}} \; dt+ \frac{2 \pi i}{p}\; z^{-\beta q/p}
\nonumber\\
& \times  \;\; \sum_{j=0}^{M-1}\,e^{(2j+1)i q \pi /p}\, \exp \!\left( -z^{-\beta/p} e^{(2j+1)i \pi/p} 
\right) \;\;, 
\label{sixtyseven}\end{align}
where it is understood that the Cauchy principal of integral is evaluated on the Stokes line, then Equivalences\ 
(\ref{sixtyfour})-(\ref{sixtysix}) can be expressed as
\begin{eqnarray}
S^{I}_{p,q}(N,z^{\beta}) \equiv F(z) + 2\pi i S z^{-\beta q/p}\, e^{-(2M+1) i q\pi/p)} \,
\exp \!\left( -z^{-\beta/p} e^{-(2M+1) i \pi/p}  \right) \;\;,
\label{sixtyeight}\end{eqnarray}
with the multiplier $S$ of the subdominant exponential term being given by
\begin{eqnarray}
S = \begin{cases} 0 \;\;,\quad -(2M+1)\pi/\beta < {\rm arg}\, z < -(2M-1) \pi/\beta \;\;, \cr
1/2 \;\;,\quad {\rm arg}\,z =-(2M+1) \pi/\beta \;\;, \cr 
1 \;\;,  \quad -(2M+3) \pi/\beta < {\rm arg}\, z< -(2M+1) \pi/\beta \;\;.
\end{cases}
\label{sixtynine}\end{eqnarray}
This is simply the conventional view as described by Berry in Ref.\ \cite{ber89} and presented at the beginning of
this section. In this instance it is clear that $S_{-}=0$ for a generalized Type I terminant when the primary 
Stokes sector is given by $|{\rm arg}\, z| < \pi/\beta$.

Now let us turn our attention to the asymptotic series $S(z)$ appearing in the complete version of Stirling's 
``approximation" for $\ln \Gamma(z)$. This series in its entirety, i.e. without any truncation whatsoever, is
found in Ref.\ \cite{kow14} to be
\begin{eqnarray}
S(z) = z \sum_{k=1}^{\infty} \frac{(-1)^k}{(2z)^{2k}}\, \Gamma(2k-1)\, c_k(1) \;,
\label{seventy}\end{eqnarray}
where the cosecant polynomials of unity \cite{kow11} can be expressed in terms of the Bernoulli numbers
as
\begin{eqnarray}
c_k(1)  = \frac{(-1)^{k}}{(2k)!} \; 2^{2k}\, B_{2k}\;.
\label{seventyone}\end{eqnarray}
Later in Ref.\ \cite{kow14} the truncation parameter $N$ is introduced with the Dirichlet series form for the
Riemann zeta function. Hence the series becomes
\begin{eqnarray}
S(z) = z \sum_{k=1}^{N-1} \frac{(-1)^k}{(2z)^{2k}}\, \Gamma(2k-1)\,c_k(1)  
- 2 z \sum_{n=1}^{\infty} S^{I}_{2,-1}\left(N,(1/2n \pi z)^2 \right)  \;.
\label{seventytwo}\end{eqnarray}
The first term on the rhs of Eq.\ (\ref{seventytwo}) is basically Paris's $\Omega(z)$, while the second term
is the frequently neglected remainder that is either divergent or conditionally convergent. Moreover, we see 
that the second term represents an infinite sum of generalized Type I terminants with $p \!=\! 2$, $q \!=\!-1$, 
$z \!=\! 1/2n\pi z$ and $\beta \!=\! 2$. By introducing these values into Equivalence\ (\ref{sixtyfour}), one obtains
\begin{align}
S^{I}_{2,-1} \left( N, (1/2n \pi z)^2 \right) &\equiv \frac{(-1)^N}{(2n\pi z)^{2N-2}} \int_0^{\infty} \frac{y^{2N-2} \,
e^{-y}}{y^2+4 n^2 \pi^2 z^2} \; dy - \frac{1}{2 n z} 
\nonumber\\
& \times \;\; \sum_{j=1}^M (-1)^{M-j}\, \exp\Bigl( -2(-1)^{M-j} n i \pi z \Bigr) \;.
\label{seventythree}\end{align}
Consequently, the asymptotic series for $\ln \Gamma(z)$ becomes
\begin{align}
S(z) &  \equiv  z \sum_{k=1}^{N-1} \frac{(-1)^k}{(2z)^{2k}}\, \Gamma(2k-1)\,c_k(1) 
 -  2  \, \Bigl( -\frac{1}{4 \pi^2 z^2} \Bigr) ^{N} \,z \sum_{n=1}^{\infty} \frac{1}{n^{2N-2}} 
\nonumber\\   
& \times \;\; \int_0^{\infty}  \frac{e^{-y}\, y^{2N-2}} {((y/2\pi z)^2 +n^2)} \; dy 
+  \sum_{n=1}^{\infty} \frac{1}{n} \sum_{j=1}^M (-1)^{M-j} \exp \left(2  (-1)^{M-j}\, n i \pi z\right) \;.
\label{seventyfour}\end{align}
Equivalences\ (\ref{seventythree}) and (\ref{seventyfour}) appear respectively as (33) and (34) in Ref.\ \cite{kow14}.
For those Stokes lines where ${\rm arg}\, (1/2n \pi z) =- (M+1/2) \pi$ or $\theta= {\rm  arg}\, z=(M+1/2) \pi$, the introduction
of Equivalence\ (\ref{sixtysix}) into Eq.\ (\ref{seventytwo}) yields 
\begin{align}
S(z)  & \equiv  z \sum_{k=1}^{N-1} \frac{(-1)^k}{(2z)^{2k}} \, \Gamma(2k-1)\,c_k(1) + 2  \Bigl( \frac{1}{4 \pi^2 |z|^2} \Bigr)^{N-1} 
z  \; \sum_{n=1}^{\infty} \frac{1}{n^{2N-2}} \; 
\nonumber\\
& \times \;\; P \int_0^{\infty} \frac{e^{-y}\, y^{2N-2}} {y^2 -4n^2 \pi^2 |z|^2} \; dy - 
i e^{i \theta} \sum_{n=1}^{\infty} \frac{1}{n} 
\sum_{j=1}^M (-1)^{j} \exp \left(2  (-1)^{j}\, n \pi |z| \right) 
\nonumber\\
& - \;\; i e^{i \theta} \sum_{n=1}^{\infty} \frac{1}{2n} \exp \left(-2 n\pi |z| \right)\;.
\label{seventyfive}\end{align}

Since each generalized Type I terminant can be expressed in terms of a Stokes multiplier accompanied by an 
exponential term that is subdominant in the vicinity of a Stokes line, it follows that $S(z)$ can be 
expressed in terms of such a multiplier and an infinite series of subdominant exponential terms as given by 
the last term on the rhs of the above equivalence. If we let 
\begin{align}
G(z) &  =  z \sum_{k=1}^{N-1} \frac{(-1)^k}{(2z)^{2k}}\, \Gamma(2k-1)\,c_k(1) 
 -  2  \, \Bigl( -\frac{1}{4 \pi^2 z^2} \Bigr) ^{N} \,z \sum_{n=1}^{\infty} \frac{1}{n^{2N-2}} 
\nonumber\\   
& \times \;\; \int_0^{\infty}  \frac{e^{-y}\, y^{2N-2}} {((y/2\pi z)^2 +n^2)} \; dy 
+  \sum_{n=1}^{\infty} \frac{1}{n} \sum_{j=1}^M (-1)^{M-j} \exp \left(2  (-1)^{M-j}\, n i \pi z\right) \;,
\label{seventyfiveb}\end{align}
where it is understood that the Cauchy principal value is evaluated at the Stokes line, then $S(z)$ can
be written as
\begin{eqnarray}
S(z) \equiv G(z) + (-1)^M S \sum_{n=1}^{\infty} \frac{e^{2 (-1)^M n i \pi z}}{n} \;\;, 
\label{seventyfivec}\end{eqnarray}
with the Stokes multiplier given by Eq.\ (\ref{sixtynine}). The above result is only partially regularized because
the series on the rhs, which arises from the infinite number of singularities situated on the Stokes line, can become 
divergent. However, the series can be regularized by using Lemma\ 2.2 in Ref.\ \cite{kow14}, which proves that 
\begin{eqnarray}
\sum_{k=1}^{\infty} \frac{ (-1)^{k+1}} {k} \; z^k \begin{cases} \equiv \ln(1+z) \quad, & \quad \Re\, z \leq -1\quad, \cr
=\ln(1+z) \quad, & \quad \Re\, z > -1 \quad . \end{cases}
\label{seventysix}\end{eqnarray} 
Then Equivalence\ (\ref{seventyfivec}) becomes
\begin{eqnarray}
S(z) \equiv G(z) + (-1)^{M+1} S \ln \Bigl( 1 - e^{2(-1)^M \pi i z} \Bigr) \;\;. 
\label{seventyseven}\end{eqnarray}
For the Stokes line where $M \!=\!0$, viz.\ $ {\rm arg}\,z \!=\! \pi/2$, the above result reduces to
\begin{eqnarray}
S(z) \equiv G(z) -  \frac{1}{2}\,  \ln \Bigl( 1 - e^{2 \pi i z} \Bigr) \;\;. 
\label{seventyeight}\end{eqnarray}
The second term on the rhs of both of the above equivalences is referred to as the Stokes discontinuity term and is denoted by 
$SD^{+}(z)$ in Eq.\ (78) of Ref.\ \cite{kow14}. It is subsequently used in the numerical study of the Stokes phenomenon in 
Sec.\ 3 of the same reference. It is also the expression which Paris claims does not represent the correct interpretation of the Stokes phenomenon. 
Yes, it is not the interpretation according to the smoothing view propounded by Berry and Olver, but it is, nevertheless, entirely 
consistent with the conventional view of the Stokes phenomenon as presented in Ref.\ \cite{kow14}.

To summarize the preceding analysis, we have seen that the subdominant exponential term in the Borel-summed regularized value of
a generalized terminant emerges as a result of a singularity being situated on a Stokes line. In the case of the asymptotic series
$S(z)$ for $\ln \Gamma(z)$, there is an infinite number of generalized terminants with coinciding Stokes lines. As a consequence, each 
Stokes line for $S(z)$ possesses an infinite number of singularities, albeit located at different positions. This produces an infinite 
number of subdominant exponential terms, each accompanied by the same Stokes multiplier that is entirely consistent with the 
conventional view of the Stokes phenomenon. The interesting property in the resulting sum of the residues due to the singularities on the
Stokes line is that it need not necessarily be convergent and thus, may require regularization. It is for this reason that the exactification 
of Stirling's approximation for $\ln \Gamma(z)$ is a challenging problem in asymptotics as discussed in the introduction of Ref.\ 
\cite{kow14}. Despite this, no one can be certain that the above analysis is indeed correct until an effective numerical study 
has been carried out. After all, we have seen that Olver claimed that his derivation of a smoothed Stokes multiplier, viz.\ ``Eq." 
(\ref{thirtyfour}), was based on rigorous mathematics. Yet he did not provide a demonstration as to just how accurate this result is. Nor
did he indicate the range of values over which it is valid. In actual fact, the notion of a proof does not really apply as far as
the Stokes phenomenon is concerned because as we shall soon see, we are talking about encountering discontinuities arising from singularities 
in Cauchy integrals and introducing an approach or method for handling them as they occur.     

If anyone possesses an incorrect interpretation of the Stokes multiplier, then it must be Paris, who seems to believe that the
Stokes multiplier should only multiply the leading subdominant exponential in a complete asymptotic expansion. Neither Stokes nor
Dingle were concerned with the leading exponential in a complete asymptotic expansion, although it must be said that they did not study 
situations where a Stokes line possesses an infinite number of singularities. In addition, the Berry/Olver derivation of Stokes smoothing presented 
earlier does not isolate the leading exponential in the subdominant part of an asymptotic solution at the expense of the other terms.
Despite this, however, it is not incorrect to isolate or factor out the leading exponential from all the other subdominant exponential terms
and then to refer to all the remaining terms as a multiplier, even though it is not what Dingle and Berry originally had in mind.  

The problem with Paris's approach like that of Berry and Olver is that in order to come up with improved results on standard
Poincar$\acute{{\rm e}}$ asymptotics, they have been forced to truncate at some stage in their analyses. As a result, we see the 
standard asymptotic constructs such as the $\sim$ symbol and $+ \cdots$, appear in (3.1) to (3.4) of Ref.\ \cite{par14}. That is, 
Paris's version of the Stokes multiplier is given by
\begin{eqnarray}
S(\theta) \sim \frac{1}{2} + \frac{1}{2} \, {\rm erf} \Bigl( c \left(\theta +\pi/2 \right) \sqrt{\pi |z|} \,\Bigr) - \frac{C_0 
e^{-2\pi \gamma |z|} }{2\pi \sqrt{z} }\;i \;\;, 
\label{seventynine}\end{eqnarray}
where $\gamma=1 + i\, \exp(i \theta)$, $C_0=B_0 \exp(-2\pi i \omega |z|) + \exp(i \omega \nu)/(1+ \exp(-i \omega))$,
\begin{eqnarray}
B_0= \frac{e^{-i \omega \alpha}}{1- e^{-i \omega}} + \frac{1}{c(\theta+\pi/2)}\;i\;\;,
\label{eighty}\end{eqnarray}
$\alpha=2 N_0 -1-2 \pi |z|$, $\nu= 2N_0-1$, $N_0$ is the optimal point of truncation and 
\begin{eqnarray}
c(\theta +\pi/2) = \omega + \frac{1}{6} \, i \omega^2 -\frac{1}{36} \, \omega^3 +\frac{1}{270} \, i \omega^4+ \cdots\;\;.
\label{eightyone}\end{eqnarray}
Consequently, the remainder is expressed as
\begin{eqnarray}
R_{N_0}(z) \approx e^{2\pi i z} \, T_{\nu} (2\pi i z)-e^{-2\pi i z} T_{\nu} (-2\pi iz)   \sim e^{2\pi i z} S(\theta) \;\;.
\label{eightyonea}\end{eqnarray}
In essence, Paris's ``refined version" of $S(\theta)$ is no different from the result given by (82) in Ref.\ \cite{kow14} except that there 
is now an extra imaginary term. This term was neglected in Ref.\ \cite{kow14} as an unnecessary complication since it is not possible to 
obtain hyperasymptotic values of the Stokes multiplier, e.g.\ to 30 figures, near the Stokes line at ${\rm arg}\,z \!=\! \pi/2$, anyway. In 
fact, the imaginary part raises another problem because imaginary parts for the multiplier simply do not appear in the conventional view of 
the Stokes phenomenon. So, if the smoothing view is correct, then it means that extra terms need to be included in the conventional view in 
order to obtain exact values of $\ln \Gamma(z)$ from its asymptotic forms. 

The whole point about the numerical study presented in Ref.\ \cite{kow14} is to determine those values of ${\rm arg}\, z$ near the Stokes
line, where Paris's form for the Stokes multiplier deviates the most from the conventional view. Specifically, these values can be
determined by plotting the real part of Approximation (\ref{seventynine}) as displayed in Figs.\ 1 of Refs.\ \cite{kow14} and \cite{par14}.  
According to Berry, Olver and Paris, one should not be able to obtain exact values of $\ln \Gamma(z)$ for these values of ${\rm arg}\,z $
via the conventional view of the Stokes phenomenon. As can be seen from the figures, the smoothed Stokes multiplier as given by 
Approximation\ (\ref{seventynine}) is closer to 1/2 than being close to 0 or 1 according to the conventional view. Hence
it stands to reason that if Stokes smoothing is correct, then the conventional view cannot possibly give accurate values of
$\ln \Gamma(z)$ in the vicinity of the Stokes line at ${\rm arg}\, z \!=\! \pi/2$.  

Table\ \ref{tab1} here or Table\ 7 in Ref.\ \cite{kow14} presents the results obtained from programming the Borel-summed regularized 
results for $\ln \Gamma(z)$ as a Mathematica module \cite{wol92}, which appears as Program 2 in the appendix of Ref.\ \cite{kow14}. 
Specifically, the module calculates $\ln \Gamma(z)$ using the following asymptotic forms:
\begin{eqnarray}
\ln \Gamma(z) = \begin{cases}  F(z) + TS_N(z) + R^{SS}_N(z) + SD^{SS,U}_{1}(z) \;,\quad   & \pi/2 < \theta \leq \pi\;, \cr
F(z) + TS_N(z) +R^{SS}_N(z) \;, \quad & -\pi/2 < \theta < \pi/2 \;. \cr
\end{cases} 
\label{eightytwo}\end{eqnarray}
In the above equation $F(z)$ represents the standard form of Stirling's approximation, i.e., 
\begin{eqnarray}
F(z) = \Bigl( z- \frac{1}{2} \Bigr) \ln z- z+\frac{1}{2} \, \ln( 2\pi)\;,
\label{eightythree}\end{eqnarray}
and the truncated series denoted by $TS_N(z)$ is given by
\begin{eqnarray}
TS_N(z) = z \sum_{k=1}^{N-1} \frac{(-1)^k}{(2z)^{2k}}\; \Gamma(2k-1)\,c_k(1) \;.
\label{eightyfour}\end{eqnarray}
As stated previously, Eq.\ (\ref{eightyfour}) is basically Paris's $\Omega(z)$. In addition, the Borel-summed remainder term denoted 
by $R^{SS}_N(z)$ was evaluated using
\begin{align}
R^{SS}_N(z) & =  \frac{2\,(-1)^{N+1} \,z}{(2\pi z)^{2N-2}} \sum_{n=1}^{\infty} \frac{1}{n^{2N-2}}\int_0^{\infty} dy \; 
\frac{y^{2N-2}\, e^{-y}}{\left( y^2 + 4 \pi^2 n^2 z^2 \right)} \;,
\label{eightyfive}\end{align}
while the final term $SD^{SS,U}_1(z)$ is the Stokes discontinuity term given by Eq.\ (\ref{seventyeight}). That is,  
\begin{eqnarray}
SD^{SS}_{1}(z) = -\ln \Bigl( 1- e^{2\pi z i} \Bigr) \;.
\label{eightysix}\end{eqnarray}

\begin{table}
\small
\centering
\begin{tabular}{|c|c|c|} \hline
$\delta$ & Method & Value \\ \hline
1/10 &   LogGamma[z]  &         -5.1085546405054331385771175 - 2.43504864133618239587613036$\, i$ \\
& $SD^{SS,U}_1(z)$ &             0.0000000146924137960847328 + 0.00000000724920978735477097$\, i$  \\
& Top &                      -5.1085546405054331385771175 - 2.43504864133618239587613036$\, i$ \\ \hline
-1/10 &  LogGamma[z]  &         -3.1156770612855851062960250 + 0.79152717486178700663566144$\, i$ \\
& Bottom &                      -3.1156770612855851062960250 + 0.79152717486178700663566144$\, i$ \\ \hline
1/100 &   LogGamma[z]  &        -4.4448078360199294879676721 - 0.68426539470619315579497619$\, i$ \\
& $SD^{SS,U}_1(z)$ &             0.0000000054543808883397577 - 0.00000000366845661861183983$\, i$ \\
& Top &                      -4.4448078360199294879676721 - 0.68426539470619315579497619$\, i$ \\ \hline
-1/100 &   LogGamma[z]  &       -4.2360547825638102221663061 - 0.35681003461125834209091866$\, i$ \\
& Bottom &                      -4.2360547825638102221663061 - 0.35681003461125834209091866$\, i$ \\ \hline
1/1000 &   LogGamma[z]  &       -4.3531757575591613140088085 - 0.53385166100905755261595669$\, i$ \\
& $SD^{SS,U}_1(z)$ &             0.0000000065016016472424544 - 0.00000000038545945628149871$\, i$ \\
& Top &                      -4.3531757575591613140088085 - 0.53385166100905755261595669$\, i$ \\ \hline
-1/1000 &   LogGamma[z]  &      -4.3322909095906129602545969 - 0.50110130347126170951651903$\, i$ \\
& Bottom &                      -4.3322909095906129602545969 - 0.50110130347126170951651903$\, i$ \\ \hline
1/10000 &   LogGamma[z]  &      -4.3438006028809735966127763 - 0.51908338527968766540121412$\, i$ \\
& $SD^{SS,U}_1(z)$ &             0.0000000065123040290213875 - 0.00000000003856476898298508$\, i$ \\
& Top &                      -4.3438006028809735966127763 - 0.51908338527968766540121412$\, i$ \\ \hline
-1/10000 &   LogGamma[z]  &     -4.3417121085407199183370966 - 0.51580834470414165478538635$\, i$ \\
& Bottom &                      -4.3417121085407199183370966 - 0.51580834470414165478538635$\, i$ \\ \hline
1/20000 &   LogGamma[z]  &      -4.3438006028809735966127763 - 0.51908338527968766540121412$\, i$ \\
& $SD^{SS,U}_1(z)$ &             0.0000000065123851251757157 - 0.00000000001928245580002624$\, i$ \\
& Top &                      -4.3438006028809735966127763 - 0.51908338527968766540121412$\, i$ \\ \hline
-1/20000 &   LogGamma[z]  &     -4.3422344065179726897501879 - 0.51662687288967352139359494$\, i$ \\
& Bottom &                      -4.3422344065179726897501879 - 0.51662687288967352139359494$\, i$ \\ \hline
\end{tabular}
\normalsize
\vspace{0.5cm}
\caption{Evaluation of $\ln\Gamma \! \left(3\exp(i(1/2+\delta)\pi ) \right)$ via Eq.\ (\ref{fiftyeight}) for
various values of $\delta$}
\label{tab1}
\end{table}

In order to be consistent with the ``Stokes smoothing" view, a relatively large value of $|z|$ was chosen when executing the module, viz.\ 
$|z| \!=\! 3$, which has an optimal point of truncation, $N_0$, approximately equal to 10. This, however, leads to another problem. 
At the present stage one does not know what form the smoothed multiplier takes for small values of $|z|$. Hence for small values of $z$, 
one cannot possibly obtain exact values of $\ln \Gamma(z)$ according to the smoothing view. Presumably, small $|z|$ implies that there is no 
optimal point of truncation, but there is no quantification from Paris on this issue. Nevertheless, no such restriction applies to any 
of the Borel-summed and MB-regularized asymptotic forms given here or in Ref.\ \cite{kow14}. That is, they are equally valid for $|z|>1$
and for $|z| \leq 1$. In addition, it should be stressed that Table\ \ref{tab1} only represents a small sample of the results from the numerical 
investigation, in which numerous values of $\delta$, where $\theta \!=\! (1/2 \!+\! \delta)\pi$, were considered. The specific values of $\delta$  
appearing in the table are those for which Stokes smoothing is expected to exhibit the greatest deviation from the step-function postulated 
in the conventional view of the Stokes phenomenon. In particular, those very close to the Stokes line are given by
$|\delta|\leq 1/100$ e.  That is, the values just below the Stokes line are given by ${\rm arg}\, z$ equal to $0.49\pi$, $0.499 \pi$, $0.4999\pi$ 
and $0.49995\pi$, while for those just above it are given by ${\rm arg}\, z$ equal to $0.51 \pi$, $0.501 \pi$, $0.5001\pi$ 
and $0.50005\pi$. Besides using different values of $\delta$ in the study, the module was also run for numerous values of 
the truncation parameter $N$. 

For each positive value of $\delta$ in Table\ \ref{tab1} there are three rows of values, while for each negative value there are 
only two rows. This is because the Stokes discontinuity term is zero for negative values of $\delta$ as indicated by Eq.\
(\ref{eightytwo}). The first row for each value of $\delta$ represents the value obtained by using the LogGamma routine in Mathematica 
and is denoted by the row with LogGamma[z] in the Method column. Depending upon whether $\delta$ is positive or not, the second row 
presents the Stokes discontinuity term according to the conventional view of the Stokes phenomenon. In general, this term was found 
to possess real and imaginary parts of the order of $10^{-8}$ or a couple of orders lower. That is, the Stokes discontinuity term is very 
small and would be either neglected or not noticeable in standard Poincar$\acute{{\rm e}}$ asymptotics, which is a direct result of choosing 
a relatively large value of $|z|$. Even though the Stokes discontinuity term is small, it is still necessary in order to give the correct 
values of $\ln \Gamma(z)$ for the hyperasymptotic calculation to thirty figures. The next value for each value of $\delta$ is labelled either 
Top or Bottom in the Method column corresponding to whether the top or bottom asymptotic form in Eq.\ (\ref{eightytwo}) has been used to calculate
$\ln \Gamma(z)$. It should also be noted that the values of the truncated sum, the regularized value of the remainder and the Stirling 
approximation were all evaluated in the Mathematica module separately, but are not displayed here due to limited space.

All the calculations carried out in the study yielded the value of $\ln \Gamma(z)$ to the hyperasymptotic accuracy as indicated in the table except  
when $|\delta|$ was extremely small, e.g.\ for $\delta \leq 10^{-5}$. Then the NIntegrate routine in Mathematica experiences convergence 
problems because the numerical integration of the remainder $R_{N}^{SS}(z)$ is too close to the singularities lying on the Stokes line. 
For example, when $\delta = 10^{-5}$, the module prints out a value of $\ln \Gamma(z)$ that agrees with the actual value to 25 decimal 
places for the real part, but in the case of the imaginary part the results only agree to 18 decimal places. The module, however, does alert
the user of the convergence problems that the NIntegrate routine experiences. Although this calculation is not presented in the table, it 
still represents a degree of success since the imaginary part of the Stokes discontinuity term is of the order of $10^{-12}$. That is, the Stokes 
discontinuity term had to be correct to the first six decimal places in order to yield the value of the imaginary part of $\ln \Gamma(z)$ 
for this very small value of $\delta$.

With the exception of the first value of $\delta$, which reflects the situation as the error function begins to veer away from
the step-function of the conventional view, we expect for all other values of $\delta$ that the real part of the smoothed Stokes 
multiplier given by Approximation (\ref{seventynine}) to be close to 1/2 according to Figs.\ 1 in Refs.\ \cite{kow14} and \cite{par14}. 
Note that Fig.\ 1 in Ref.\ \cite{par14} gives the Stokes multiplier for $|z|= 8$, which has an optimal point of truncation that is 
approximately equal to 26. As a consequence, the transition is far more rapid than in Fig.\ 1 of Ref.\ \cite{kow14}. Whilst Paris presents  
values of the Stokes multiplier for ${\rm arg}\, z $ ranging form 0.325 to 0.750 at intervals of 0.025, he does not give the values for the
Stokes multiplier in the vicinity of the Stokes line as in Table\ \ref{tab1}. Nevertheless, for ${\rm arg}\, z= 0.475$, he obtains
a value of $0.2894310-0.0182669\,i$ for his smoothed multiplier, while for ${\rm arg}\, z= 0.525$, he obtains a value of
$0.7105689-0.0182669\,i$. Although these are not representative of the situation in Table\ \ref{tab1} and moreover, are nowhere near
the accuracy needed to conduct the hyperasymptotic investigation in Table\ \ref{tab1}, they indicate that the Stokes
discontinuity term should be at least 30 percent less than the figures in the table for $\delta >0$ and that thirty percent of 
the Stokes discontinuity needs to be added to the values in which $\delta < 0$ according to the smoothing view of the Stokes phenomenon.   
That is, the top asymptotic form in Eq.\ (\ref{eightytwo}) with about half the Stokes discontinuity term should be a far more accurate 
approximation to the actual value of $\ln \Gamma(z)$. However, we see the opposite. The first asymptotic form yields the exact value of 
$\ln \Gamma(z)$ for all values of $\delta$ greater than zero despite the fact that the Stokes discontinuity term has no effect on the 
first nine decimal places. For $\delta < 0$, according to the smoothing view the bottom asymptotic form in Eq.\ (\ref{eightytwo}) should 
not yield exact values of $\ln \Gamma(z)$ because it is missing about half the Stokes discontinuity term. Once again, we observe the 
opposite; the bottom asymptotic form yields exact values of $\ln \Gamma(z)$ for all negative values of $\delta$ in the table. Thus, 
it is evident that there is no smoothing of the Stokes phenomenon occurring in the vicinity of the Stokes line at $\theta = \pi/2$ or
else it would not been possible to give the exact values of $\ln \Gamma(z)$ from Eq.\ (\ref{eightytwo}). 

From the results in the table we see that by using the conventional view of the Stokes multiplier we are able to obtain thirty
figure accuracy for $\ln \Gamma(z)$ in the vicinity of the Stokes line. Moreover, we could have considered more decimal places, if
this was really necessary. To achieve higher levels of accuracy all one has to do is alter the working precision plus 
the precision and accuracy goals in the Mathematica module. This will come at the expense of the computing time. Despite this,
has Paris done the same in Ref.\ \cite{par14}? If the smoothing point of view for the Stokes phenomenon is indeed 
correct, then it should not only be vastly superior to the conventional view by yielding more accurate values of
$\ln \Gamma(z)$, it should also expose where the errors or deficiencies occur in the conventional view. After all, there is simply 
no point in offering an alternative view to the mathematical community if it is unable to provide an improvement on the existing 
view/approach. Despite this, there is no numerical evidence in Ref.\ \cite{par14} demonstrating how smoothing
is able to match or even provide more accurate values of $\ln \Gamma(z)$ in the vicinity of a Stokes line. All we see 
is a table giving the Stokes multiplier (both real and imaginary parts) to six decimal figures for a much larger value
($|z|=8$) than in Table\ \ref{tab1}  because he was unable to obtain this level of accuracy for the Stokes multiplier
for $|z|\!=\!3$.  In addition, the analysis resulting in Table\ \ref{tab1} was conducted for numerous values of the truncation parameter 
$N$ far away from the optimal point whereas Paris's computations require an optimal point of truncation. This means that his approach is
useless for $|z|<1$, a limitation that should never apply in (hyper)asymptotics. 

The issues mentioned in the preceding paragraph are the ones that need to be addressed by Paris. Until his smoothing approach is able 
to yield hyperasymptotic values of a special function irrespective of whether the variable is large or small and without the need 
for optimal truncation, there is no point adopting the concept of smoothing to asymptotics, which is not only vague and limited, but 
also very unwieldy. At the very least Paris needs to match the results of Table\ \ref{tab1} with his smoothed forms using
$|z|=3$, not $|z|=8$ as in Ref. \cite{par14}, which possesses a much higher optimal point of truncation. 

\section{Conclusion}
Although the Stokes discontinuity term presented here would be neglected in standard Poincar$\acute{\rm e}$ asymptotics, it raises the 
question of how can a function that is continuous develop asymptotic forms that are discontinuous. Doubtless to say, 
Stokes \cite{sto04} was the first to be intrigued by this behaviour. This puzzling behaviour arises because on a Stokes line asymptotic 
series for $\ln \Gamma(z)$ are maximally divergent \cite{din73}. In another words, all its terms are real and of the same sign. For such 
values of ${\rm arg}\, z$, divergent series are particularly counter-intuitive. For example, we have seen that the regularized value of the 
geometric series is given by $1/(1-z)$. For positive real values where $z>1$, the geometric series is positive definite, but its regularized 
value is negative. 

We have seen that the discontinuities along a Stokes line arose from the Borel summation of asymptotic series along a Stokes line. 
This regularization technique produces Cauchy integrals which yield different values according to whether their singularities lie above, on 
or below the line contour along the positive real axis. In order to obtain exact values of the original function via its asymptotic forms, one 
requires the regularized remainder of all the asymptotic series despite the fact that for large values of $|z|$ truncation may yield a good 
approximation. Even then the relative accuracy of the truncated series varies for each value of $z$ and $N$. Thus, across a Stokes line Borel 
summation of an asymptotic series results in a discontinuous jump in the multiplier accompanying the subdominant term. For $|z|> 1$ 
we have seen that these discontinuous terms are tiny or subdominant and are neglected in standard Poincar$\acute{\rm e}$ asymptotics. Nevertheless, 
they are important in a hyperasymptotic calculation. So, in response to the above question, it is because we are dealing with the regularization 
via Borel summation of divergent series. Consequently, we should not expect divergent series to behave in a manner that is consistent with 
convergent series. To make matters more mysterious, the emergence of jump discontinuities when the original function is continuous does not occur 
when asymptotic series undergo Mellin-Barnes regularization as described in Ref.\ \cite{kow14}.  Yet, the resulting asymptotic forms from
both regularization techniques yield identical values of $\ln \Gamma(z)$ for all values of $z$ including arguments. 

In summary, the two main claims in Ref.\ \cite{par14} that (1) there is no need for the concept of regularization and (2)
the ``Stokes smoothing" view originally due to Berry \cite{ber89} represents the correct view for interpreting asymptotic behaviour 
near a Stokes line, have been refuted. In regard to the first claim, we have seen that unbeknownst to him, Paris has actually employed 
the concept of regularization to arrive at the forms that he derives for the remainder $R_{N}(z)$ given by Eqs.\ (\ref{four})-(\ref{fiveb}) 
here. On the second claim, we have seen that the Borel-summed regularized values of $\ln \Gamma(z)$ in Ref.\ \cite{kow14} are based 
entirely on the conventional view of the Stokes phenomenon, which holds that the Stokes multiplier of the subdominant part of a complete 
asymptotic expansion behaves as a step-function at a Stokes line. Moreover, Paris has not provided a hyperasymptotic study demonstrating 
that his smoothing view is indeed superior to the conventional view. Specifically, he needs to reproduce or better still improve upon the 
results presented in Table\ \ref{tab1}, which display the results from a hyperasymptotic analysis of $\ln \Gamma(z)$ in Ref.\ \cite{kow14} 
based on the conventional view of the Stokes phenomenon for $|z|=3$. After all, if the smoothing approach is indeed correct, then this should 
become evident when one wishes to obtain hyperasymptotic values of a function from its asymptotic forms not only near, but also far away from 
a Stokes line, regardless of whether optimal truncation has been invoked or not. It should also be able to handle all values of $|z|$, not 
just $|z| \gg 1$.

\end{document}